\documentclass[11pt]{article}

\usepackage[top= 2cm,bottom=2cm,left=2.00cm,right=2.00cm]{geometry}
\usepackage{amsfonts,amsmath,graphicx,color}
\usepackage{mathrsfs}
\makeatletter

\renewcommand{\d}{{\rm d}}
\newcommand{\e}{{\rm e}}
\renewcommand{\i}{{\rm i}}
\newcommand{\PD}[2]{\frac{\partial #1}{\partial #2}}
\newcommand{\FD}[2]{\frac{\d #1}{\d #2}}

\newcommand{\D}{{\rm D}}

\DeclareMathSymbol{\ZSet}{\mathalpha}{AMSb}{"5A}
\DeclareMathSymbol{\RSet}{\mathalpha}{AMSb}{"52}
\DeclareMathSymbol{\CSet}{\mathalpha}{AMSb}{"43}

\title{Synchrony in firing rate neural networks with multiple delays: A harmonic balance approach}
\author{S. Coombes\footnote{School of Mathematical Sciences, University of Nottingham, Nottingham.
NG7 2RD, UK. \newline {\tt email:} stephen.coombes@nottingham.ac.uk} \and H.G.E. Meijer\footnote{Department of Applied Mathematics, University of Twente, PO Box 217, 7500 AE Enschede, The Netherlands. \newline {\tt email:} h.g.e.meijer@utwente.nl}}
\date{}

\begin{document} 
\maketitle

\begin{abstract}
\noindent
Networks of neural mass nodes with delayed interactions are increasingly being used as models for large-scale brain activity. To complement the growing number of computational studies of such networks, it is timely to develop new mathematical studies of their solution structure and bifurcations. The analysis of steady states and their stability is relatively well developed, though that of time-periodic solutions is far less so. Here, we show how the method of harmonic balance is ideally suited to describing delay-induced and delay-modulated periodic oscillations at both the node and network level. This approach reduces the formally infinite dimensional setting of the delayed differential equation network to a finite dimensional one, opening the way for a practical combined analytical and numerical treatment. At the node level, we show how to construct periodic orbits and develop an associated linear stability analysis to determine the Floquet exponents and, thus, stability. At the network level, we further show that explicit progress for analysing the stability of the synchronous state can be made for networks with a circulant structure. This is achieved with the use of an adjacency lag operator that decouples the linearised network equations into a set of equations that can each be analysed using the techniques previously developed for analysing a single node. When combined with numerical continuation techniques this allows us to build the skeleton of a network bifurcation diagram, and highlight the role of distance-dependent delays in contributing to novel spatio-temporal patterns arising from the instability of a synchronous state, including travelling periodic waves, alternating anti-phase solutions (in which only next nearest neighbours are synchronised), cluster states, and more exotic behaviours. The predictions from this harmonic balance framework are found to be in excellent agreement with direct simulations.
\end{abstract}

\section{Introduction}

The use of neuroimaging modalities such as electro- and magneto-encephalography has been used to great effect in cognitive psychology and medicine to shed light on brain dynamics in both functional and dysfunctional contexts \cite{Gross2019}. In particular, functional connectivity reflecting temporal correlations between spatially separated brain regions can be related to behaviour in a variety of different tasks, and is known to have a neural basis. This has prompted several modelling studies to help understand the role of local neuronal population dynamics and network structure in contributing to large-scale brain dynamics. The common ingredients for such studies are a low-dimensional neural mass model for the local node dynamics, a set of structural interactions, and a corresponding set of delays. The neural mass model is used as a proxy for the behaviour of a large population of spiking excitatory and inhibitory cells, albeit idealised using coarse-grained variables (say for population averages) and sacrifices the details of spike times in favour of firing rates. A classic example is the Wilson--Cowan model consisting of two coupled ordinary differential equations \cite{Wilson72}. The use of diffusion magnetic resonance imaging now gives us an invaluable source of information about the structural organisation of the brain. This includes both the strength of connectivity between anatomically defined brain regions and data about axonal projections that mediate long-range interactions between distant neuronal ensembles. The combination of neural mass models and this connectome data forms the basis for several recent studies, as exemplified by the work of Deco \textit{et al}. \cite{Deco2009}, Abeysuriya \textit{et al}. \cite{Abeysuriya2018}, and Castaldo \textit{et al}. \cite{Castaldo2023}.
The latter emphasises that the conduction delays arising from the finite propagation speed of signals along axonal pathways are vitally important for capturing realistic functional connectivity patterns and their time evolution, and builds on previous work that already suggested that coupled oscillators with delayed interactions could give rise to transient moments of synchronisation between clusters of nodes in a network at frequencies that are lower than the intrinsic frequency of oscillation of uncoupled nodes \cite{Cabral2022}. Indeed, many computational studies of model oscillators (not necessarily of neural mass type) have illustrated that transmission time delays can play a major role in organising network behaviour, such as synchronisation and its disruption \cite{Kutchko2013,Petkoski2019,Conti2019,Nakagawa2023}.

In contrast to the number of computational studies mentioned above, there is relatively little corresponding mathematical analysis beyond that of steady state behaviour. This can perhaps be attributed to the challenge of even analysing a dynamical system with just one delay\footnote{The solution space for delay differential equations has to be considered as infinite-dimensional although only a finite number of dynamical variables is involved. As a consequence, nonlinear delay differential equations reveal a broad class of instabilities leading from oscillatory to chaotic behaviour, and see e.g., \cite{Erneux2009} in an applied mathematics context.} before even worrying about how to address network states. Regarding steady state behaviour, many results are now known, relying mainly on linear stability analysis, and see \cite{Campbell2007} for an excellent overview. Moreover, for single nodes (with a small number of dynamical variables) it is practical to use numerical tools such as DDE-BIFTOOL \cite{Engelborghs2002} to track the stability of equilibrium points and find bifurcations to periodic orbits, as done in \cite{Coombes2009} for a single Wilson--Cowan node with two delays. A more recent study of a Wilson--Cowan network has shown that delay-induced instabilities of the network steady state can excite patterns consistent with human resting state functional connectivity (in the $\alpha$ frequency band) \cite{Tewarie2019}, and see also \cite[Ch. 8]{Coombes2023}. The topic of this paper is, however, the analysis of \textit{oscillatory} network states. One natural way to analyse these is to use the theory of weakly coupled oscillators \cite{Hoppensteadt97}, with a reduction to a phase description as in \cite{Ton2014} or assuming a phase-oscillator network of Kuramoto type as in \cite{Petkoski2016}. However, this is valid only for weak interactions and has to be augmented in the case that the oscillation at the node level is delay induced\footnote{A numerical technique for computing the phase response function for delay differential equations can be found in \cite{Kotani2012,Novicenko2012} and for the computation of phase \textit{and} amplitude response see \cite{Nicks2024}.}.
Here, we shall consider strongly coupled networks of Wilson--Cowan nodes with multiple time-delayed interactions. For tractability, we will restrict attention to networks with distance-dependent interactions.
Recent work has shown how techniques from nonsmooth dynamics can be used to handle the special case that the firing rate function in the Wilson--Cowan model is a Heaviside \cite{Sayli2024}. Here, we develop a more general approach that is valid for smooth firing rate functions.

In Section \ref{Sec:model} we introduce the delayed differential equations representing the neural network model. Next, in Section \ref{Sec:harmonicbalance} we review the method of harmonic balance and show how it can be used for the construction of periodic orbits in systems with delay. Moreover, we show how this constructive approach can naturally feed in to a Floquet stability analysis, leading to the determination of Floquet exponents as zeros of a certain complex function. By way of illustration we demonstrate the utility of this approach for the construction and stability of periodic orbits in a scalar delay differential equation (DDE) with a single delay. We then bootstrap this to the case of a Wilson--Cowan node (with a 2-dimensional state vector). Numerical continuation codes, as employed in MatCont \cite{MATCONT}, are used in conjunction with the harmonic balance method to build bifurcation diagrams and show the dependence of solution branches under variation of the delay.
Section \ref{Sec:Ring} treats the main topic of this paper, namely the stability of the synchronous state in a network where the strength of connections and the delays between nodes are distance-dependent.
This is achieved by building on Section \ref{Sec:harmonicbalance} with a suitable `diagonalisation' of the network Floquet theory, allowing for the construction of a skeleton bifurcation diagram organised around solution branches of the synchronous network state. Predictions from this bifurcation analysis are shown to be in excellent agreement with direct numerical simulation of the neural network model\footnote{Direct simulations of the delayed network model were performed in Julia using the {\tt DifferentialEquations} package \cite{Rackauckas2017}.}.
Finally, in Section \ref{Sec:Discussion} we revisit some of the main results obtained and discuss extensions for further work.

\section{Model equations\label{Sec:model}}

A single Wilson--Cowan node is built from an interacting set of dynamic excitatory and inhibitory activities. Each of these activities is described by a scalar so that the node has a 2-dimensional state vector. The variable $u$ will be used to denote the activity of the excitatory population and $v$ that of the inhibitory population. The relative time-sale between the dynamics for each activity will be denoted by $\kappa$.
The nonlinear interaction between the excitatory and inhibitory populations is mediated by a nonlinear, typically sigmoidal, firing rate function $F$. With the introduction of an integer index $i=1,\ldots,N$ then a network of interacting Wilson--Cowan nodes can be written as
\begin{equation}
\begin{aligned}
 \label{network}
\FD{}{t} u_i(t) &= -u_i(t) + F \left (I_u+ w^{uu}u_i(t-\tau_0)-w^{vu}v_i(t-\tau_0) + \epsilon \sum_{j=1}^N w_{ij} u_j(t-\tau_{ij}) \right ), \\
\kappa \FD{}{t} v_i(t) &= - v_i(t) + F \left (I_v + w^{uv} u_i(t-\tau_0) - w^{vv} v_i(t-\tau_0) \right) , \qquad i=1,\ldots,N.
\end{aligned}
\end{equation}
Here, we have assumed that all nodes in the network are identical, namely that there is no network index dependence of $\kappa$ and $F$.
The parameters $w^{\alpha \beta} \geq 0$, $\alpha,\beta \in \{ u,v\}$, specify the \textit{local} connection strengths and $\tau_0 \geq 0$ specifies a \textit{self-delay}.
The structural interactions between nodes $i$ and $j$ is captured by the connection strength $w_{ij} \geq 0$ and the corresponding delay by $\tau_{ij} \geq 0$. The connections between different nodes are only through their excitatory components, consistent with the large scale wiring pattern of the brain. The parameter $\epsilon \geq 0$ sets the global scale for the strength of interaction. Throughout the rest of this paper we work with the choice
\begin{equation}
F(x) = \frac{1}{1+\e^{-\beta x}}, \qquad \beta >0 .
\label{Activation}
\end{equation}
The constants $I_\alpha \in \RSet$ may be thought of as background drives or as effective thresholds.
A single node (recovered from (\ref{network}) with the choice of a fixed index $i$ and $\epsilon=0$) can support a periodic orbit via a Hopf bifurcation in the absence of a self-delay ($\tau_0=0$). Before destabilising the steady-state via a Hopf bifurcation it is also possible to generate a delay induced oscillation with a non-zero self-delay. For more about the bifurcation structure of the Wilson--Cowan node model in the absence and presence of delays see \cite{Coombes2009}. If the network has a graph Laplacian structure ($w_{ij} = -a_{ij} + \delta_{ij} \sum_k a_{ik}$ for some arbitrary set of connections $a_{ij} \in \RSet$) or a row-sum constraint 
($\sum_{j} w_{ij} = \text{const}$ for all $i$) then (\ref{network}) admits synchronous solutions of the form $(u_i,v_i) = (\overline{u},\overline{v})$ for all $i=1,\ldots,N$, which could be both time-independent and time-dependent.
Regarding the network analysis of (\ref{network}) this has mainly been limited to the linear stability analysis of a time-independent homogeneous steady state \cite{Tewarie2019}, the study of the master stability function for determining the stability of the synchronous time-periodic network state in the absence of delays \cite{Ahmadizadeh2016,Coombes2018} and an extension of this approach for a single small delay \cite{AlDarabsah2021}, all for networks with either a row-sum constraint or Laplacian coupling. A recent numerical bifurcation analysis (using DDE-BIFTOOL) of a small ring of Wilson--Cowan oscillators with a single delay can be found in \cite{Pinder2023}. Here the authors find that for weak coupling transitions between synchronous and alternating anti-phase states (in which only next nearest neighbours are synchronised) occur via a torus bifurcation (with a route to chaos via a torus breakdown), and that with stronger coupling multi-stable dynamics can occur.
In the following sections we shall develop a bifurcation analysis of synchronous time-periodic states that can be delay-induced or delay-modulated for the network model (\ref{network}) with multiple time-delayed interactions. To guarantee the existence of a synchronous solution we shall focus on networks with a circulant structure in both the strength of connections and their delays. Namely, we shall work with the choice $w_{ij} = w_{|i-j|}$ and $\tau_{ij} = \tau_{|i-j|}$.

Although we use the Wilson--Cowan model as an exemplar of a neural mass the analysis that we present here is agnostic to the choice of node dynamics and could equally well be applied to other neural mass models, such as that in \cite{Byrne2022}, or indeed any DDE model with smooth nonlinearities. Moreover, the methodology can also be applied to the analysis of periodic travelling waves that arise in networks with non-delayed interactions, as recently considered by Ruschel and Giraldo \cite{Ruschel2025}.

\section{Analysis: harmonic balance\label{Sec:harmonicbalance}}

The harmonic balance method has been widely used for analysing the periodic solutions of nonlinear differential equations, often in a mechanical setting, see e.g., \cite{Detroux2015}. Roughly speaking a truncated Fourier series representation is used to convert a set of differential equations into a nonlinear algebraic system of equations for the (complex) Fourier coefficients. These can then be found using numerical algorithms such as the Gauss--Newton method. It is also sometimes named as the Fourier--Galerkin method, since it consists in the application of the Galerkin method with Fourier basis and test functions. In essence, periodic functions are approximated with their Fourier coefficients, which become the new unknowns of the problem. The term harmonic balance was first introduced by Krylov and Bogoliubov \cite{Krylov1943}, and importantly the method can also be used to study periodic solutions to delay equations, see e.g., \cite{MacDonald1995,Liu2010,Sun2023}. Here, we use this method as a starting point for developing an analytical approach for the Floquet theory of the network model (\ref{network}). Other starting points could be based upon alternative methods for approximating periodic solutions to dynamical systems such as the Poincar\'e--Lindstedt or Shohat methods as in \cite{Simmendinger1999}, though we will argue that the choice of harmonic balance is both elegant and pragmatic.
To aid in the exposition of this new methodology for treating periodic solutions in networks with multiple delays we shall first give a \text{primer} that focuses on a scalar DDE with a single delay. We then move on to a single Wilson--Cowan node before treating the full network in Section \ref{Sec:Ring}.

\subsection{A primer\label{Sec:Primer}}

Let us first consider an idealisation of the Wilson--Cowan model that treats only one of the excitatory or inhibitory populations with a self-feedback term. This gives rise to scalar DDE with a single delay of the form
\begin{equation}
\FD{}{t} x(t) = -x + f(x(t-\tau)), \qquad x \in \RSet, \qquad \tau \geq 0 ,
\label{xprimer}
\end{equation}
where $f(x) = F(I+wx)$ and $w, I \in \RSet$.
A $T$-periodic orbit $x_0(t)=x_0(t+T)$ can be represented with a Fourier series as:
\begin{equation}
x_0(t) = \sum_{n \in \ZSet} a_n \e^{\i \omega_n t}, \quad \omega_n = \frac{2 \pi n}{T}, \quad a_n \in \CSet .
\label{Fourierx}
\end{equation}
Since $x_0(t) \in \RSet$ then $a_{-n}=a_n^*$, where $*$ denotes complex conjugation.
If we restrict the sum in (\ref{Fourierx}) to $n \in \{-M, \ldots, 0, \ldots, M \}$ then there are $2M+2$ unknowns to determine, namely $a_0 \in \RSet$, $a_{1,\ldots,M} \in \CSet$, and $T$.
To determine these we use the harmonic balance method and sample $x_0(t)$ at the $2M+1$ instants, $t_n=nT/(2M+1)$, $n=-M,\ldots,0,\ldots,M$:
\begin{equation}
x_0(t_n) = \sum_{m=-M}^M a_m \e^{2 \pi \i n m/(2M+1)} .
\end{equation}
Introducing the notation $X=(x_0(t_{-M}), \ldots, x_0(t_0), \ldots, x_0(t_M))$ and $A=(a_{-M}, \ldots, a_0, \ldots, a_M)$, allows us to write the relation between the sampled states and the Fourier coefficients as the linear system:
\begin{equation}
X = S A ,
\label{X=SA}
\end{equation}
where $S$ is the Vandermonde matrix:
\begin{equation}
[\, S \,]_{nm} = \exp(2 \pi \i n m/(2M+1)).
\label{S}
\end{equation}
From (\ref{xprimer}), the periodic orbit $x_0(t)$ satisfies
\begin{equation}
\left ( 1 +\FD{}{t} \right ) x_0(t) = f(x_0(t-\tau)).
\label{x0}
\end{equation}
Using (\ref{Fourierx}) the left hand side of (\ref{x0}) takes the form $\sum_{m \in \ZSet} \left ( 1 +\i \omega_m \right ) a_m \e^{i \omega_m t}$.
Sampling at times $t_n$ and truncating the sum leads to a nonlinear systems of $2M+1$ algebraic equations for the components of $X$:
\begin{equation}
S \Omega S^{-1} X - f(S \Gamma S^{-1} X) = 0,
\label{Xa}
\end{equation}
where 
\begin{equation}
[\, \Omega \,]_{nm} = \delta_{nm} \left ( 1 +\i \omega_m \right ) , \qquad [\, \Gamma \,]_{nm} = \delta_{nm} \e^{-\i \omega_m \tau} .
\end{equation}
Here, $S^{-1}$ can be calculated explicitly as $S^{-1} = S^*/(2M+1)$. To fix an origin of the one parameter family of periodic orbits we set the condition $x_0'(0)=0$, or equivalently the single equation
\begin{equation}
(-M,\ldots,0,\ldots,M)^{\mathsf{T}} S^{-1} X = 0 .
\label{Xb}
\end{equation}
The $2M+2$ set of equations given by (\ref{Xa}) and (\ref{Xb}) determine the $2M+2$ unknowns, and may be solved numerically (e.g., with {\tt fsolve} in Matlab). The complex Fourier coefficients are then constructed from the real vector $X$ as $A=S^{-1}X$.

To determine the stability of the periodic orbit $x_0(t)$, consider a small $T$-periodic perturbation $z(t)$ such that $x(t) = x_0(t) + \e^{\lambda t}z(t)$, where $\lambda \in \CSet$.
This perturbation satisfies the linearised equation:
\begin{equation}
\left ( \lambda + 1 + \FD{}{t} \right ) z(t) = \e^{-\lambda \tau} f'(x_0(t-\tau)) z(t-\tau) .
\label{zt}
\end{equation}
Differentiation of (\ref{x0}) with respect to $t$ and comparison with (\ref{zt}) shows that $x_0'(t)$ is the eigenfunction of the linearised system with $\lambda=0$.
Sampling (\ref{zt}) at times $t_n$ and using a truncated Fourier series representation of $z(t)$ as $z(t) = \sum_{m=-M}^M z_m \e^{\i \omega_m t}$ yields
\begin{equation}
\left [ S \mathcal{L}(\lambda) S^{-1} - \e^{-\lambda \tau} \D f(S \Gamma S^{-1}X) S \Gamma S^{-1} \right ] Z =0 ,
\label{Z}
\end{equation}
where $\D f (X) =\operatorname{diag} (f'(X))$, $Z=(z(t_{-M}), \ldots, z(t_0), \ldots, z(t_M))$, 
$F'(x) = \beta F(x) (1-F(x))$, and
 \begin{equation}
 [\, \mathcal{L}(\lambda) \,]_{nm}= \delta_{nm} \left (\lambda+1 +\i \omega_m \right ) .
\end{equation}
For non-trivial solutions of (\ref{Z}) we require that $\mathcal{E}(\lambda)=0$, where 
\begin{equation}
\mathcal{E}(\lambda) = \det (S \mathcal{L}(\lambda) S^{-1} - \e^{-\lambda \tau} \D f(S \Gamma S^{-1}X) S \Gamma S^{-1} ).
\end{equation}
The periodic orbit $x_0$ is stable if $\text{Re} \, \lambda <0$. Note that the method presented here is only an approximation of the orbit and its spectrum, with errors expected to reduce as $M \rightarrow \infty$.

\begin{figure}[htbp]
\begin{center}
\includegraphics[width=.75\paperwidth]{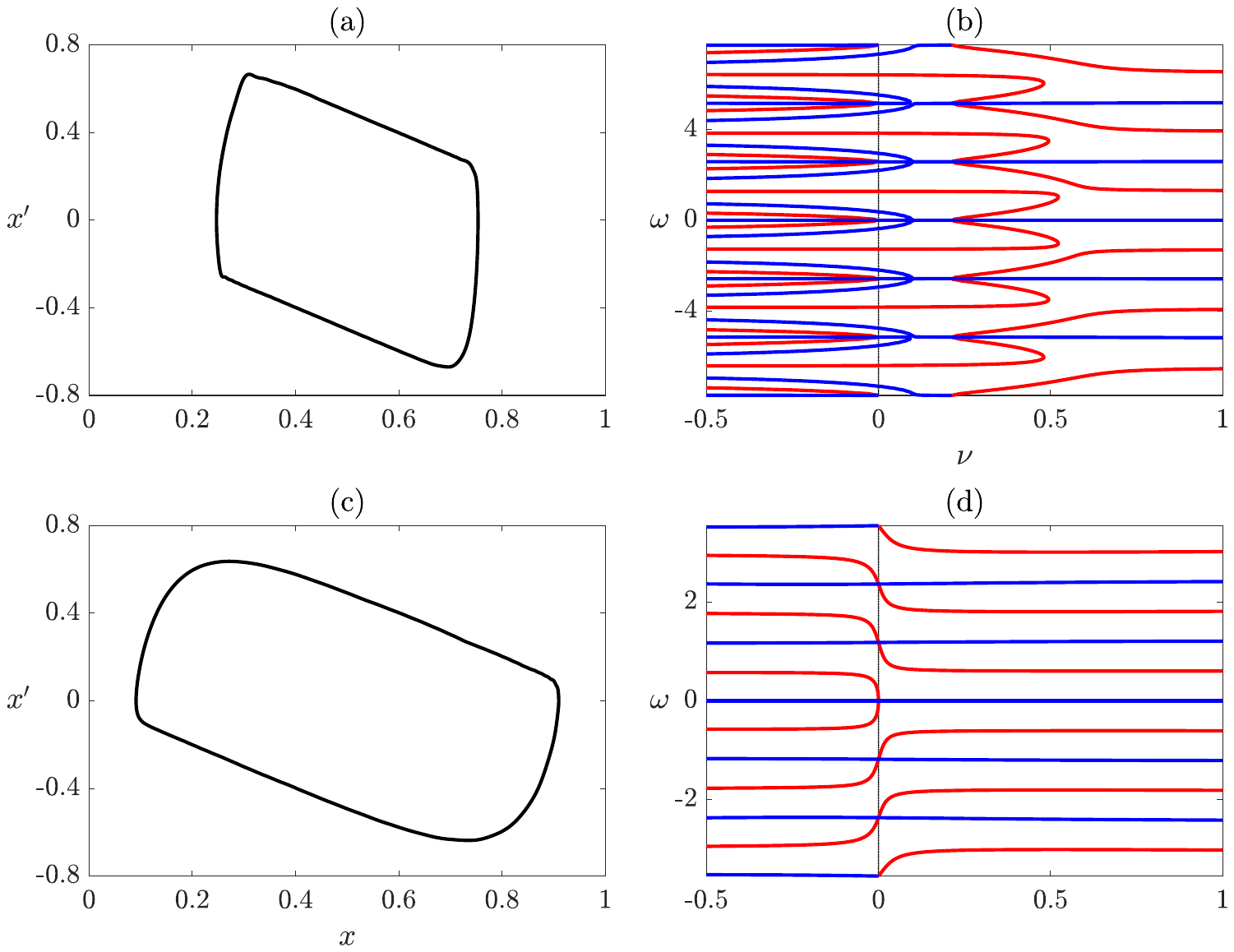}
\caption{Construction and stability of periodic orbits for the \textit{primer} model given by (\ref{xprimer}) with $\tau=2$ and $M=50$. Left: Periodic obits in the $(x,x')$ plane. Right: Stability plots in the complex $(\nu,\omega)$ plane with eigenvalues determined by the crossing of the red and blue lines where $\text{Re} \, \mathcal{E}(\nu + \i \omega) = 0$ and $\text{Im} \, \mathcal{E}(\nu + \i \omega) = 0$ respectively. (a) An unstable periodic orbit as predicted from the spectrum in panel (b), where eigenvalues can be found in the right hand complex plane. Parameters: $\beta=80$, $w=1$, and $I=-1/2$.
(c) A stable periodic orbit as predicted from the spectrum in panel (d), with no eigenvalues in the right hand complex plane. Parameters: $\beta=20$, $w=-1$, and $I=1/2$. The zero eigenvalue seen in both panels (c) and (d) reflects time-translation invariance of the periodic orbits.}
\label{Fig:Primer}
\end{center}
\end{figure}

An example of the use of this methodology to construct periodic orbits and determine their stability is shown in Fig.~\ref{Fig:Primer}. 


\subsection{Wilson--Cowan node with a single delay}
\label{WCsingledelay}

For a single Wilson--Cowan node the technique described above in Section \ref{Sec:Primer} generalises naturally with the use of tensor notation to cope with the increased dimensionality of the local dynamics. To illustrate this more concretely consider a system with a state vector $x \in \RSet^p$ with a single delay $\tau_0$ such that the dynamics can be written as the DDE:
\begin{equation}
\FD{}{t} x(t) = G (x(t), x(t-\tau_0)) .
\label{G}
\end{equation}
Consider a $T$-periodic orbit $x_0(t) = x_0(t+T)$ with a Fourier series representation that generalise (\ref{Fourierx}) with the choice $a_n \in \CSet^p$.
Using the notation $X=(x_0(t_{-M}), \ldots, x_0(t_0), \ldots, x_0(t_M)) \in \RSet^{(2M+1)p}$ and $A=(a_{-M}, \ldots, a_0, \ldots, a_M) \in \CSet^{(2M+1)p}$, allows us to write the relation between the sampled states and the Fourier coefficients as the linear system (\ref{X=SA}) under the replacement
\begin{equation}
S \rightarrow S_p := S \otimes I_p ,
\end{equation}
where $\otimes$ represents the Kronecker tensor product and $I_p$ is the $p \times p$ identity matrix.
Following the process described in the primer, for sampling and truncation, leads to a nonlinear systems of $(2M+1)p$ algebraic equations for the components of $X$:
\begin{equation}
S_p L_p S_p^{-1} X - G(X,S_p \Gamma_p S_p^{-1} X) = 0,
\label{Xaa}
\end{equation}
where $L_p:=L \otimes I_p$ and $\Gamma_p:=\Gamma \otimes I_p$, and
\begin{equation}
[\, L \,]_{nm} = \delta_{nm} \left ( \i \omega_m \right ) , \qquad [\, \Gamma \,]_{nm} = \delta_{nm} \e^{-\i \omega_m \tau_0} .
\end{equation}

To fix an origin of the one parameter family of periodic orbits we set the $n$th component of $x$ to have a vanishing derivative at the origin, or equivalently: 
\begin{equation}
e_n^{\mathsf{T}} \left \{ (-M,\ldots,0,\ldots,M)^{\mathsf{T}} S_p^{-1} X \right \} = 0 ,
\label{origin}
\end{equation}
where $e_n \in \RSet^p$ is a canonical vector for the $n$th direction.
The determination of stability also proceeds as in the primer and it is convenient to introduce the Jacobians $J_j(t) \in \RSet^{p \times p}$ as:
\begin{equation}
J_j(t) = \left. \PD{}{x_j} G(x_1,x_2) \right |_{(x_1,x_2) = (x_0(t), x_0(t-\tau_0))}, \qquad j=1,2 .
\end{equation}
The periodic orbit $x_0$ is stable if $\text{Re} \, \lambda <0$, where $\mathcal{E}(\lambda)=0$, and
\begin{equation}
\mathcal{E}(\lambda) = \det (S_p \mathcal{L}_p(\lambda) S_p^{-1} - \D J_1 - \e^{-\lambda \tau_0} \D J_2 \, S_p \Gamma_p S_p^{-1} ).
\end{equation}
Here, $\mathcal{L}_p(\lambda) = \mathcal{L}_p(\lambda) \otimes I_p$, $[\, \mathcal{L}(\lambda) \, ]_{nm}= \delta_{nm} \left (\lambda+\i \omega_m \right )$ and $\D J _j = \operatorname{block diag} (J_j(t_{-M}), \ldots, J_j(t_{0}), \ldots, J_j(t_{M})).$

For a single Wilson--Cowan node we have that $p=2$ and
\begin{equation}
G((u,v), (u_{\tau_0},v_{\tau_0}) ) = \begin{bmatrix} -u + F \left (I_u+ w^{uu}u_{\tau_0}-w^{vu}v_{\tau_0} \right ) \\ \left (- v + F \left (I_v + w^{uv} u_{\tau_0} - w^{vv} v_{\tau_0} \right ) \right )/ \kappa \end{bmatrix} ,
\end{equation}
where $x_\tau$ is shorthand for $x(t-\tau)$. This gives $J_1(t) =- \Upsilon$, where $\Upsilon$ is a $2 \times 2$ time-independent matrix given by 
\begin{equation}
\label{Upsilon}
\Upsilon = \begin{bmatrix} 1 & 0 \\0 & 1/\kappa \end{bmatrix}, 
\end{equation}
and
\begin{equation}
J_2(t) = \begin{bmatrix}
w^{uu} F' \left (I_u+ w^{uu}u(t-\tau_0)-w^{vu}v(t-\tau_0) \right ) & -w^{vu} F' \left (I_u+ w^{uu}u(t-\tau_0)-w^{vu}v(t-\tau_0) \right ) \\
w^{uv} F' \left (I_v + w^{uv} u(t-\tau_0) - w^{vv} v(t-\tau_0) \right )/\kappa & - w^{vv} F' \left (I_v + w^{uv} u(t-\tau_0) - w^{vv} v(t-\tau_0) \right ) /\kappa
\end{bmatrix} .
\end{equation}

\begin{figure}[htbp]
\begin{center}
\includegraphics[width=.75\paperwidth]{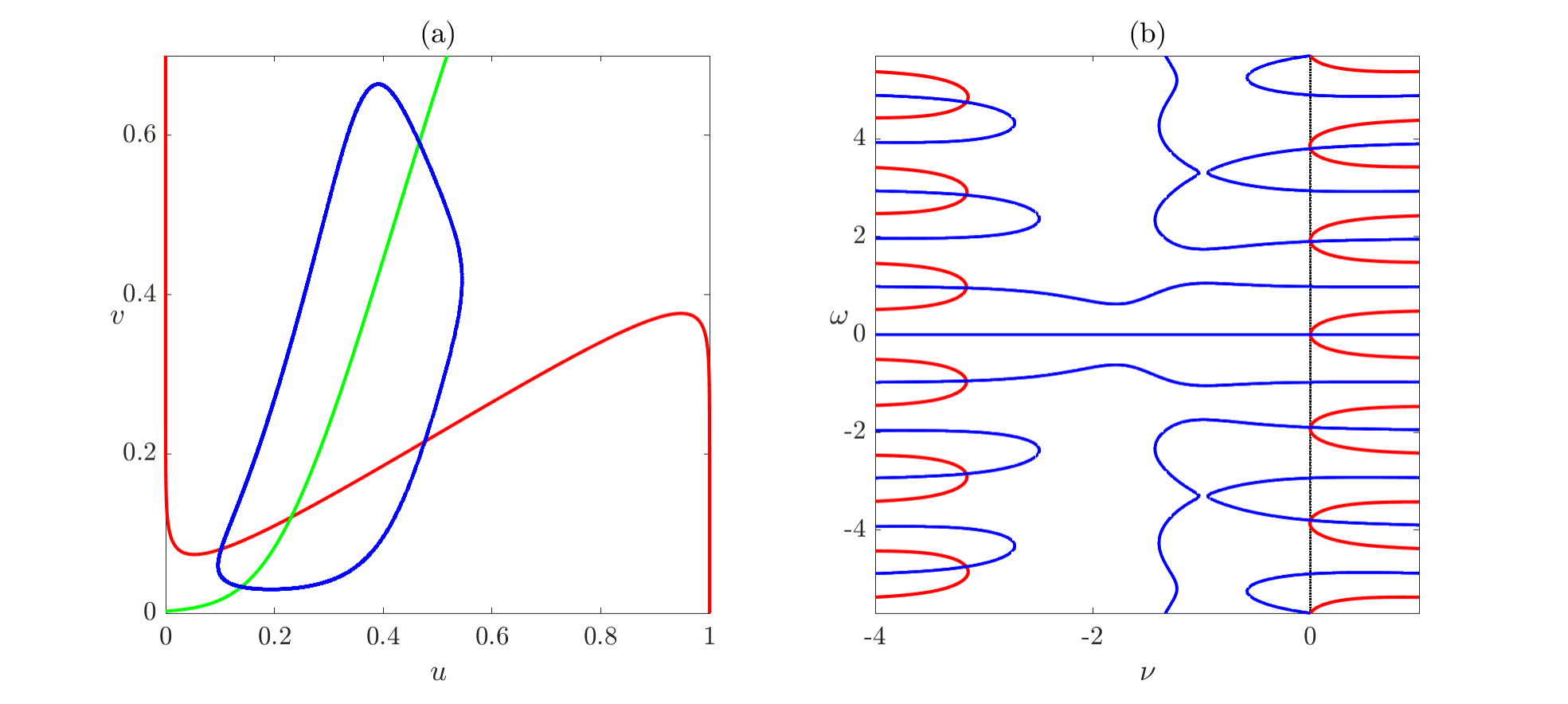}
\caption{A stable delay induced periodic orbit in a single Wilson--Cowan node. 
(a) A periodic orbit constructed using the harmonic balance method with $M=30$ is shown in blue.
The red and green curves denote the $v$ and $u$ nullcline respectively for the model in the absence of delay.
(b). A plot of the spectrum for the orbit shown in (a), with eigenvalues determined by the crossing of the red and blue lines where $\text{Re} \, \mathcal{E}(\nu + \i \omega) = 0$ and $\text{Im} \, \mathcal{E}(\nu + \i \omega) = 0$ respectively. Since there are no eigenvalues in the right hand complex plane the orbit is linearly stable.
Parameters: $\tau_0 = 0.2$, $\kappa=1/2$, $w^{uu}= 1$, $w^{vu} = 2$, $w_{vv} = 1/4$, $w^{uv}= 1$, $I_u=-0.05$, $I_v=-0.3$, $\beta=20$.
\label{Fig:WilsonCowan}
}
\end{center}
\end{figure}

An example of a periodic orbit for a single Wilson--Cowan node constructed using the machinery above is shown in Fig.~\ref{Fig:WilsonCowan}. Here, the orbit is purely delay induced and is linearly stable since no eigenvalues are found in the right hand complex plane. In the absence of delay ($\tau_0=0$), and for other parameters as fixed in Fig.~\ref{Fig:WilsonCowan} (see caption), a standard linear stability analysis of the (single) steady state shows that it is a stable focus.

%
%

\section{Ring networks and an adjacency lag operator\label{Sec:Ring}}

Here we shall consider ring networks of $N$ nodes with distance between nodes $i$ and $j$ defined by $\operatorname{dist}(i,j)=\min(|i-j|, N-|i-j|)$.
Moreover, we will assume that that network interactions can be described by circulant matrices such that 
\begin{align}
w_{ij} &= w(\operatorname{dist} (i,j) ), \\
\tau_{ij} &= 
	\begin{cases}
		\operatorname{dist} (i,j)/c & \operatorname{dist} (i,j) > 0 \\
		\tau_0 & \operatorname{dist} (i,j) = 0
	\end{cases} ,
\end{align}
where $c$ represents a fixed and uniform action potential conduction speed.
Substituting a synchronous solution of the form $(u_i(t), v_i(t)) = (\overline{u}(t), \overline{v}(t))$, for all $i=1, \dots, N$, into the network model (\ref{network}), shows that, if it exists, then a synchronous network state is a periodic solution of 
\begin{equation}
\begin{aligned}
	 \label{sync}
	\FD{}{t} \overline{u} (t) &= -\overline{u}(t)+ F \left (I_u+ w^{uu}\overline{u}(t-\tau_{0}) -w^{vu} \overline{v}(t-\tau_{0}) + \epsilon \sum_{j=0}^{N-1} w_j \overline{u}(t-\tau_{j}) \right ), \\
	\kappa \FD{}{t} \overline{v}(t) &= -\overline{v}(t)+ F \left (I_v+ w^{uv} \overline{u}(t-\tau_0) - w^{vv} \overline{v}(t-\tau_0) \right) ,\end{aligned}
\end{equation}
where $w_j = w( \operatorname{dist} (j,0))$ for all $j=0,\ldots,N-1$ and $\tau_j= \operatorname{dist}(j,0)/c$ for $j = 1,\ldots,N-1$.

The harmonic balance method may once again be used to make progress in the construction and stability of the synchronous network state.
By extending the results in section \ref{WCsingledelay}, the period $T$ and periodic orbit $x_0(t) = (\overline{u}(t), \overline{v}(t))$ can be found as the solution to (\ref{origin}) and
\begin{equation}\label{Eq:RingCont}
S_p L_p S_p^{-1} X - G(X,S_p \Gamma_p (0) S_p^{-1} X, S_p \Gamma_p (1) S_p^{-1} X, \ldots, S_p \Gamma_p (N-1) S_p^{-1} X) = 0,
\end{equation}
where $G(x(t),x(t-\tau_0),x(t-\tau_1), \ldots, x(t-\tau_{N-1}))$ is a generalisation of the right hand side of (\ref{G}) to account for a multitude of delays with $\Gamma_p(k) = \Gamma(k)\otimes I_p$ (with $p=2$ for the Wilson--Cowan model) and
\begin{equation}
[\, \Gamma(k) \,]_{nm} = \delta_{nm} \e^{-\i \omega_m \tau_k}. 
\end{equation}

In order to analyse stability of synchronous orbit we consider perturbed solutions in the form $(u_i(t),v_i(t))=(\overline{u}(t)+\delta u_i(t), \overline{v}(t)+\delta v_i(t))$. Then perturbations evolve according to 
\begin{align}
	\FD{}{t}\delta u_i(t) &= -\delta u_i(t) +F'(\chi_{\overline{u}}) \left ( w^{uu} \delta u_i(t-\tau_0) - w^{vu} \delta v_i(t-\tau_0) +\epsilon \sum_{j=1}^N w_{ij} \delta u_j(t-\tau_{ij}) \right ), \label{deltaui} \\
	\kappa \FD{}{t}\delta v_i(t) &= - \delta v_i (t)+ F'(\chi_{\overline{v}}) \left (w^{uv} \delta u_i(t-\tau_0) - w^{vv} \delta v_i (t-\tau_0) \right),
	\label{deltavi}
\end{align}
where $\chi_{\overline{u}}(t) = I_u+ w^{uu}\overline{u}(t-\tau_{0}) -w^{vu} \overline{v}(t-\tau_{0}) + \epsilon \sum_{j=0}^{N-1} w_j \overline{u}(t-\tau_{j})$ and $\chi_{\overline{v}}(t) = I_v+ w^{uv} \overline{u}(t-\tau_0) - w^{vv} \overline{v}(t-\tau_0)$.
We introduce a $2 \times 2$ block notation using:
\begin{align}
	W_0(t) =\begin{bmatrix}   F' (\chi_{\overline{u}}(t)) (w^{uu} + \epsilon w_0) &  -F' (\chi_{\overline{u}}(t)) w^{vu}  \\ 
	F' (\chi_{\overline{v}}(t)) w^{uv} &  -F' (\chi_{\overline{v}}(t))w^{vv}  \end{bmatrix},
	\quad
	W_k(t) = \epsilon \begin{bmatrix}   F' (\chi_{\overline{u}}(t)) w_k &  0  \\ 
	0 &  0  \end{bmatrix}, \quad k=1, \dots, N-1.
\end{align}
We can then rewrite (\ref{deltaui})-(\ref{deltavi}) in the vector form for $i=1, \dots, N$ as: 
\begin{align}
&	\FD{}{t}\begin{bmatrix} \delta u_1(t) \\ 	\kappa \delta v_1(t) \\ \vdots \\  \delta u_N(t) \\  	\kappa \delta v_N(t)  \end{bmatrix}=-\begin{bmatrix} \delta u_1(t) \\ 	 \delta v_1(t) \\ \vdots \\ \delta u_N(t) \\  	\delta v_N(t)  \end{bmatrix}+
	\begin{bmatrix} {W_0}(t) &0&0&\dots &0 \\ 
		0& {W_0}(t) &0&\dots &0 \\ 
		\vdots& \dots&\ddots&\dots &\vdots \\ 
		0&0&\dots & {W_0}(t) &0 \\ 
		0&\dots &0&0&{W_0}(t) \end{bmatrix}\begin{bmatrix} \delta u_1(t-\tau_{0}) \\  	 \delta v_1(t-\tau_{0}) \\ \vdots \\  \delta u_N(t-\tau_{0}) \\  	\delta v_N(t-\tau_{0})  \end{bmatrix} \nonumber \\ &+
		\begin{bmatrix} 0&{W_1}(t) &0&\dots & 0  \\  
		0 & 0&{W_1}(t) &\dots &0  \\  
		\dots& \dots&\dots&\dots &\dots \\ 
		0&\dots&0 &0&{W_1}(t)  \\  
		{W_1}(t) &\dots &0&0 &0 \end{bmatrix}\begin{bmatrix} \delta u_1(t-\tau_{1}) \\  	 \delta v_1(t-\tau_{1}) \\ \vdots \\  \delta u_N(t-\tau_{1}) \\  	\delta v_N(t-\tau_{1})  \end{bmatrix}+\dots \nonumber \\&+	\begin{bmatrix} 0&\dots&0&0 &0 &0&0&\dots &{W_{N-1}}(t)   \\  
		{W_{N-1}}(t) &\dots&0&0&0 &0 &0&\dots & 0 \\ 
		\dots&{W_{N-1}}(t)&\dots&\dots&\dots&\dots&\dots& \dots &\dots\\  
		0 & 0 &0&\dots&\dots&\dots&\dots& \dots &0\\  
		\dots&\dots&\dots&\dots&\dots&\dots&{W_{N-1}}(t)& \dots &\dots\\  
		0&\dots& 0 & 0 &0&\dots&0&  {W_{N-1}}(t)&0\\   \end{bmatrix}\begin{bmatrix} \delta u_1(t-\tau_{N-1}) \\  	 \delta v_1(t-\tau_{N-1}) \\ \vdots \\  \delta u_N(t-\tau_{N-1}) \\  	\delta v_N(t-\tau_{N-1})  \end{bmatrix}. \label{ComponentNotation}
\end{align}
We note that each of the square matrices on the right hand side of (\ref{ComponentNotation}) has a circulant structure that can be written in terms of a generator matrix $C$ in the form 
$\operatorname{blockdiag}(W_k(t), \ldots, W_k(t)) (C^k \otimes I_p)$, where $C \in \RSet^{N \times N}$ is given by
\begin{equation}
C = \begin{bmatrix}
0 & I_{N-1} \\
1 & 0
\end{bmatrix} .
\end{equation}
We now use a vector notation with $Y=(u_1, v_1, u_2, v_2, \ldots , u_N, v_N) \in \RSet^{2N}$ and denote the synchronous solution by $\overline{Y}(t)= (\overline{u}(t), \overline{v}(t), \overline{u}(t), \overline{v}(t), \ldots, \overline{u}(t), \overline{v}(t))$ and perturbations $\delta Y = (\delta u_1, \delta v_1, \delta u_2, \delta v_2, \ldots , \delta u_N, \delta v_N)$.
Then, (\ref{ComponentNotation}) can be written succinctly in the vector notation
\begin{equation}
	\FD{}{t} \delta Y(t) = - (I_N \otimes \Upsilon) \delta Y(t) + \left [\sum_{k=0}^{N-1}(C^k \otimes \Upsilon W_k (t) ) \, \mathcal{A}_k \right ]\delta Y(t) ,
	\label{deltaY}
\end{equation}
where $\Upsilon$ is given by (\ref{Upsilon}) and $\mathcal{A}_k$ is the \textit{lag operator}: $\mathcal{A}_k u(t) = u(t-\tau_k)$. It is thus natural to identify $C^k \mathcal{A}_k$ as an \textit{adjacency lag operator} (since it treats lags at various distances $k$ across the network) \cite{Otto2018}.

To `diagonalise' (\ref{deltaY}) we make use of the fact that all the circulant matrices $C^k$ (which effectively encode the topology of the ring network) have a common eigenspace.
The eigenvectors of a circulant matrix have the form $e_q = (1,\omega_q, \omega_q^2, \ldots, \omega_q^{N-1})/\sqrt{N}$, where $q=0,\ldots, N-1$, and $\omega_q=\exp(2 \pi \i q /N)$ are the $N$th roots of unity. 
The eigenvalues of $C$ are given explicitly by $\omega_q$ and the eigenvalues of $C^k$ are thus $(\omega_q)^k$

If we introduce the matrix of eigenvectors $P=[e_0 ~e_1 ~\ldots ~e_{N-1}]$, then we have that 
\begin{equation}
	P^{-1} C^k P=\operatorname{diag}(\omega_0^k, \omega_1^k, \ldots, \omega_{N-1}^k) \equiv \Lambda_k.
\end{equation}
Using a change of variable $\delta Z=(P\otimes I_p)^{-1}\delta Y$ and applying $(P\otimes I_p)^{-1}$ to both sides of the system (\ref{deltaY}) we obtain
\begin{equation}\label{variationalnetwork2}
		\FD{}{t}\delta Z (t) = -(I_N\otimes \Upsilon) \delta Z (t) 
		+ \left [\sum_{k=0}^{N-1}(\Lambda_k \otimes \Upsilon W_k (t) ) \, \mathcal{A}_k \right ]\delta Z(t) . 
\end{equation}
This has a $p \times p$, $N-$block structure with the dynamics in each block given by 
\begin{equation}\label{BlockManyD}
	\FD{}{t}\xi_q(t)= \Upsilon\left[-I_m + \sum_{k=0}^{N-1} W_k (t) \omega_q^k \, \mathcal{A}_k \right] \xi_q(t), \quad q=0,\dots, N-1,\quad \xi_{q} \in \mathbb{C}^2.
\end{equation}
Thus, the synchronous orbit will be linearly stable provided all solutions of $\mathcal{E}_q(\lambda)=0$ have $\operatorname{Re} \lambda<0$ for all $q=0,\ldots,N-1$ (excluding a parameter independent zero eigenvalue arising from time-translation invariance), where
\begin{equation}\label{Eq:RingStab}
\mathcal{E}_q(\lambda) = \det \left (S_p \mathcal{L}_p(\lambda) S_p^{-1} - \D J_{-1} - \sum_{k=0}^{N-1} \D J_{k} \, S_p \Gamma_p(k) S_p^{-1} \, \omega_q^k \e^{-\lambda \tau_k} \right ).
\end{equation}
Here, $\mathcal{L}_p(\lambda)$ is given in Section \ref{WCsingledelay} (with $p=2$ for the Wilson--Cowan model) and 
\begin{align}
\D J_{-1} & = \operatorname{blockdiag} (- \Upsilon, \ldots, -\Upsilon, \ldots, -\Upsilon), \\
\D J_k & = \operatorname{blockdiag} (\Upsilon W_k(t_{-N}), \ldots, \Upsilon W_k(t_{0}), \ldots, \Upsilon W_k(t_N)). \label{DJ}
\end{align}

%

\subsection{Numerical results}
Next, we use our mathematical framework in conjunction with numerical tools to illustrate the utility of our approach in uncovering some of the spatio-temporal patterning that can arise in delayed networks with distance-dependent coupling.  For illustrative purposes, we restrict attention to the case of $N=7$ nodes, though our methodology applies to ring networks of arbitrary size. We predict linear instabilities for certain modes and use direct simulations of the nonlinear model to validate these and uncover the full emergent behaviour. We start with a network version of the primer model, and proceed to the Wilson--Cowan case making full use of the tensor approach.   For notational clarity we introduce two new symbols for better distinguishing between the self-delay and the interaction delays.  Namely, we introduce  $\tau_{\text{intra}}:= \tau_0$ and $\tau_{\text{inter}} = \tau_{j} - \tau_{j-1}$ for all $j=2, \ldots,N-1$.


\subsubsection{An excitatory ring network\label{numerics_primer}}
We choose a purely excitatory network (obtained from (\ref{network}) by excluding all inhibitory variables $v_i$) with connectivity decaying with distance according to $\epsilon w_{ij}=A C g^{-\operatorname{dist}(i,j)}$, with $g=0.5$.  We choose $C$ such that $\sum_{j}w_{ij}=A$ to compare for different numbers of nodes. 
For the activation function we choose $\beta=10$ in (\ref{Activation}) and for the self-delay we fix $\tau_{\text{intra}}=1.5$.
We vary the inter-node delay $\tau_{\text{inter}}$ to explore the effect of distant-dependent coupling on synchrony. For the harmonic balance method, we choose $M=30$ modes. While the branches of synchronous and travelling wave solutions are already well represented for lower values of $M$, their stability is more accurately described for $M=30$. And results for $M=40$ were similar. Combining the harmonic balance method with the phase condition to have a maximum at $t=0$ defines a regular system of equations that can be solved for the coefficients $X$ and period $T$. Adding a pseudo-arclength condition allows numerical continuation in a system parameter \cite{MeijerDercoleOldeman2009}. We use the secant method to approximate the tangent vector. After having computed the branches, we determined the stability at each point. Initially, we visually inspect contour lines of the real and imaginary parts of the determinant (\ref{Eq:RingStab}). The matrix elements are well-defined, but the determinant takes large values. To monitor stability along a branch of synchronous solutions, we compute its zeros using bisection on the real axis and as contour intersections for complex roots.

\begin{figure}[htbp]
\begin{center}
\includegraphics[width=.75\paperwidth]{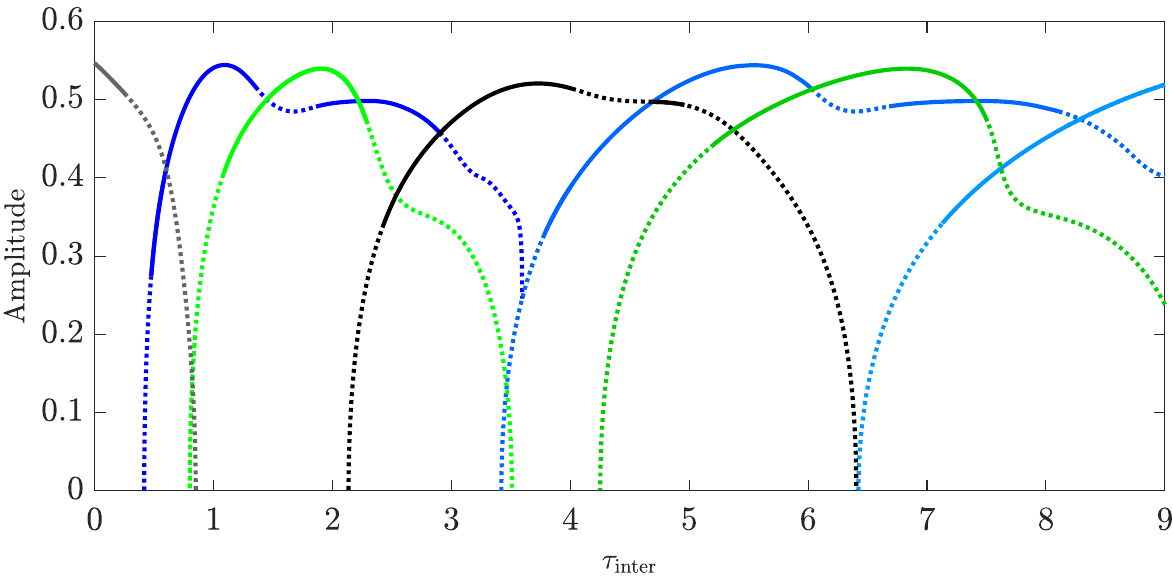}
\caption{Branches of periodic solutions in an excitatory  ring network computed with numerical continuation. Blue: Synchronous solution, Green: Travelling wave solution with a $q=1$ mode pattern, Black: Travelling wave solution with a $q=2$ mode pattern, and grey for $q=3$ near $\tau_{\rm inter}=0.2$. Solid parts of the branches are stable, while the dotted parts indicate unstable oscillations. Because of the delay, any periodic solution is again a periodic solution if the period $T$ is added to the delay, i.e. $\tau_{\text{inter}}\rightarrow T+\tau_{\text{inter}}$. This fact explains the similar shape of the secondary branches.}
\label{Fig:PrimerBif}
\end{center}
\end{figure}

\begin{figure}[ht!]
\begin{center}
\includegraphics[width=16cm]{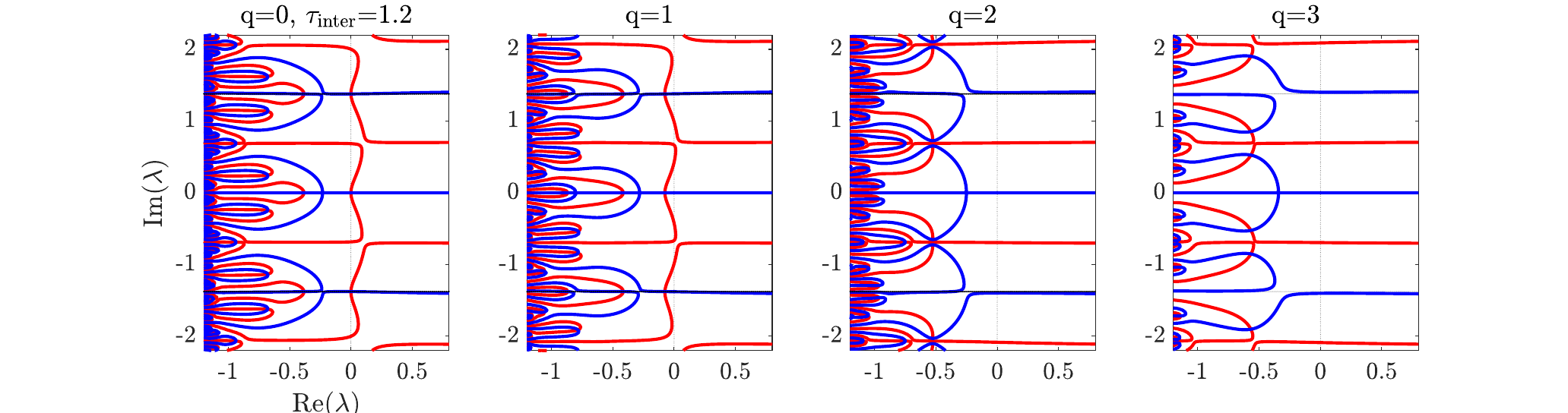}\\
\includegraphics[width=16cm]{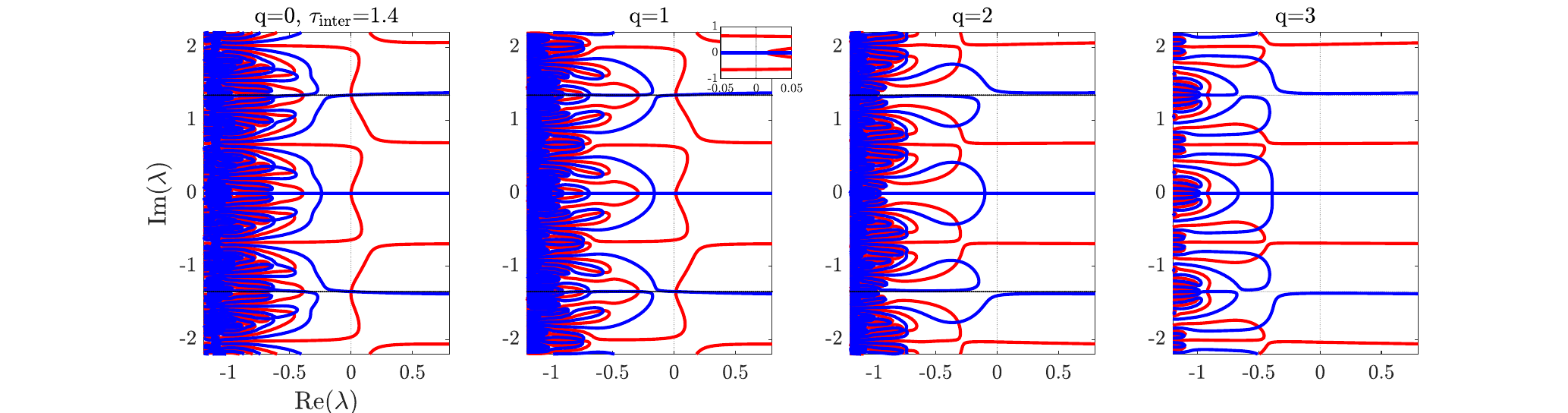}
\caption{
(Top) $\tau_{\text{inter}}=1.2$; all intersections of the red and blue contour lines (zero level sets of the real and imaginary parts of $\mathcal{E}_q(\lambda)$ given by (\ref{Eq:RingStab})) are in the left half plane. 
(Bottom) $\tau_{\text{inter}}=1.4$; one intersection for $q=1$ of the contour lines is in the right half plane. The inset shows $\lambda=0.0154$ is an unstable eigenvalue.}
\label{fig:N7spectrum}
\end{center}
\end{figure}

\begin{figure}[ht!]
\includegraphics[width=7cm]{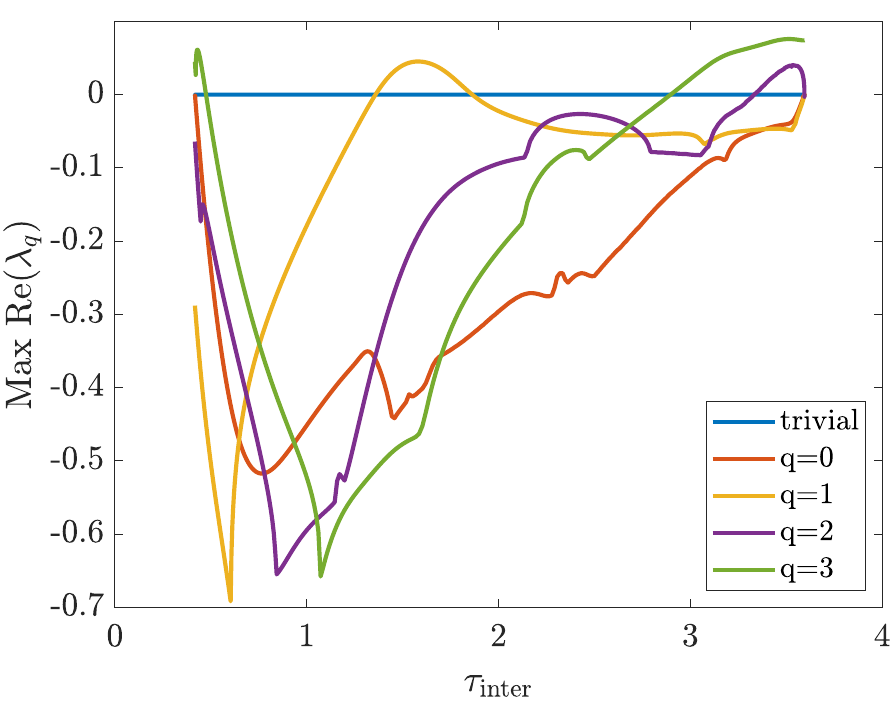}
\includegraphics[width=10cm]{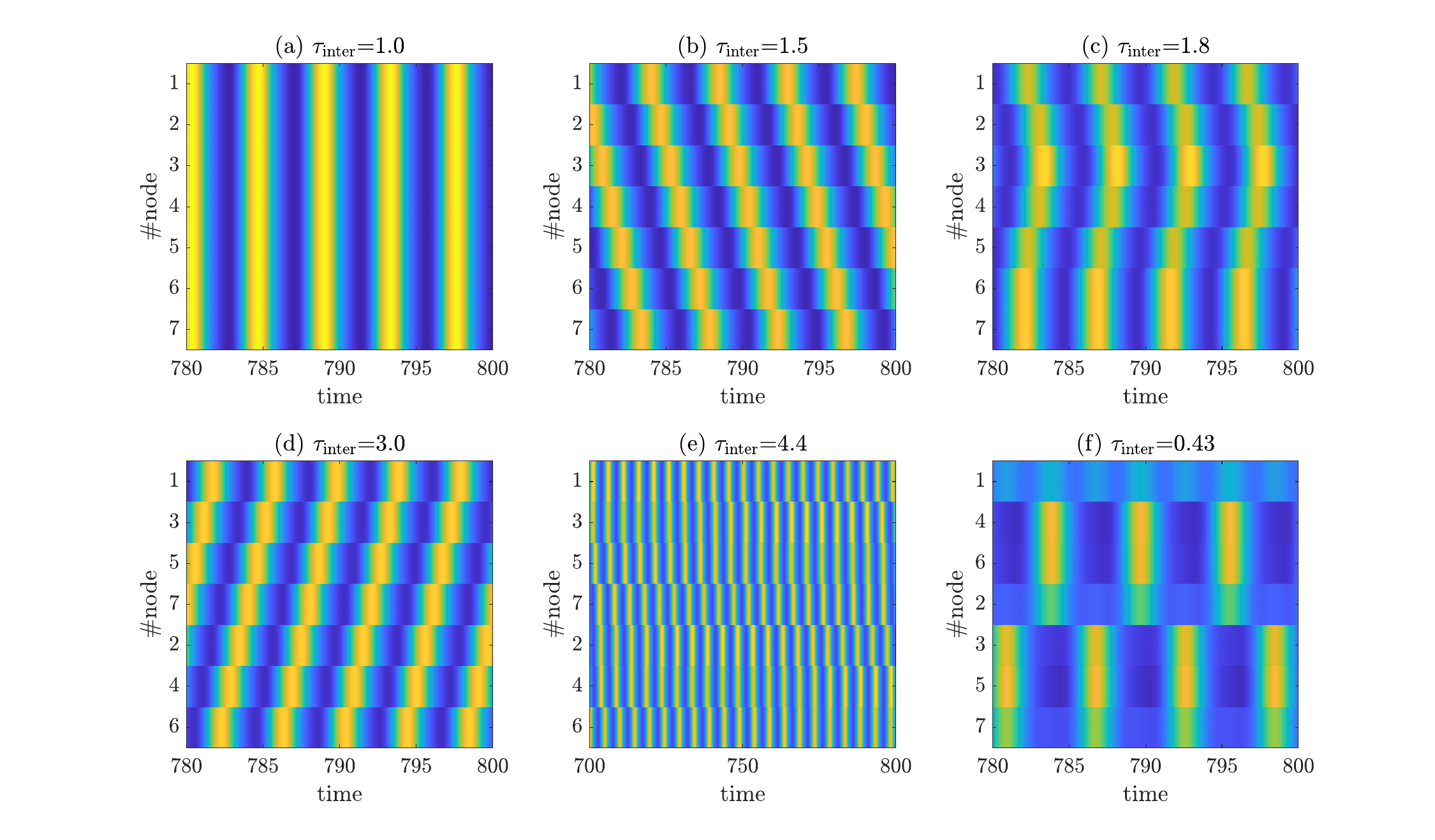}
\caption{(Left) For each $q$, the eigenvalue(s) with largest real part, including the trivial one for $q=0$ along the entire primary synchronous branch are shown.
(Right) Stable patterns from direct simulations: (a) synchronous solution with $\tau_{\text{inter}}=1.0$, (b) travelling wave ($q=1$) with $\tau_{\text{inter}}=1.5$, (c) slightly asynchronous oscillation breaking up into two waves with $\tau_{\text{inter}}=1.8$, (d) travelling wave $q=2$ where we reordered the nodes on the vertical axis to highlight the wave with $\tau_{\text{inter}}=2.8$, (e) travelling wave with disordered phase with $\tau_{\text{inter}}=4.6$, (f) out-of-phase oscillation breaking up into two clusters with $\tau_{\text{inter}}=0.43$.}
\label{fig:N7patterns}
\end{figure}

In Fig.~\ref{Fig:PrimerBif}, we show a bifurcation diagram obtained using this setup. We first computed the branches using (\ref{Eq:RingCont}), starting at $\tau_{\text{inter}}=1.0$, where we find a fully synchronous solution (Fig.~\ref{fig:N7patterns}a). Using numerical continuation, we follow this solution for increasing and decreasing $\tau_{\text{inter}}$. This branch starts from a Hopf bifurcation at $\tau_{\text{inter}}\approx 0.4176$ and ends at another synchronous branch at $\tau_{\text{inter}}\approx 3.5956$. Along this branch, we show how an instability arises. At $\tau_{\text{inter}}=1.2$ and $\tau_{\text{inter}}=1.4$, we have evaluated (\ref{Eq:RingStab}) and plotted the zero contours of the real and imaginary parts of
$\mathcal{E}_q(\nu+\i \omega)$ 
for $q=0,1,2,3$, see Fig.~\ref{fig:N7spectrum}. For $\tau_{\text{inter}}=1.2$, we see that $\lambda=0$ is a solution for $q=0$ as expected from translation invariance. All other intersections of the red and blue contours occur in the left half plane, indicating all exponents are negative showing  that the synchronous solution is stable. When we increase $\tau_{\text{inter}}=1.4$, we see that an instability arises for $q=1$, see inset zooming in near the origin, with $\lambda=0.0154$. A simulation with a small perturbation shows that this $q=1$ mode grows slowly, but eventually the pattern shifts from synchronous to a travelling wave.

We have systematically computed the spectrum for the synchronous solution. For each $q$, we determined the eigenvalue with the largest real part, see Fig.~\ref{fig:N7patterns}. Increasing $\tau_{\text{inter}}$, we find an instability at $\tau_{\text{inter}}=1.36$ for $q=1$. If we add a small perturbation with that mode, this then grows until the system settles to a travelling wave (and Fig.~\ref{fig:N7patterns}b predicts that perturbations with $q=2,3$ return to the synchronous solution). This travelling wave is a separated branch and not a nearby solution, as we have traced with numerical continuation. Here we used the same harmonic balance equation (\ref{Xa}) and (\ref{Xb}) but with an additional phase shift for the input from the other nodes (on top of the one due to the delay). Using simulations, we determined this branch is stable between $\tau_{\text{inter}}=1.06$ and $\tau_{\text{inter}}=2.29$. There is a different scenario when we follow the synchronous solution decreasing from $\tau_{\text{inter}}=2.0$. At $\tau_{\text{inter}}\approx 1.87$, there is an instability of mode $q=1$ as confirmed with our stability analysis, and see Fig.~\ref{fig:N7patterns} (Left). Initially, the seven nodes oscillate only slightly out-of-phase with a similar amplitude as seen in Fig.~\ref{fig:N7patterns}c. Decreasing $\tau_{\text{inter}}$ further, two clusters seem to develop, each with some nodes having a large amplitude, and others less so. Then at $\tau_{\text{inter}}\approx 1.65$, this pattern becomes unstable and the travelling wave appears again. For much larger values of the inter-node delay $\tau_{\text{inter}}\in[2.41,2.90]\cup[3.80,4.96]$, the synchronous solution co-exists with another travelling wave (mode $q=2$, i.e. next-nearest neighbour pattern). This mode has an additional phase modulation of the oscillation travelling over the network for $\tau_{\text{inter}}\in[4.02,4.68]$, see Fig.~\ref{fig:N7patterns}e. Increasing $\tau_{\text{inter}}$ further, we encounter similar scenarios as in many other systems with a delay $\tau$, namely a periodic solution with period $T$ is again a solution for delay $\tau+T$.

Restarting at $\tau_{\text{inter}}=1.0$ and decreasing this inter-node delay, we see first that the amplitude of the oscillation decreases. Then at $\tau_{\text{inter}}\approx 0.475$, there is an instability with $q=3$, and here the network splits into two clusters but with various amplitudes, 
see Fig.~ \ref{fig:N7patterns}f. At the transition from these clusters to a travelling wave, we observe for $\tau_{\text{inter}}\approx 0.265$ a travelling modulation on top of the clustered oscillation.

\subsection{Wilson--Cowan ring network} 
\begin{figure}[ht!]
\includegraphics[width=6cm]{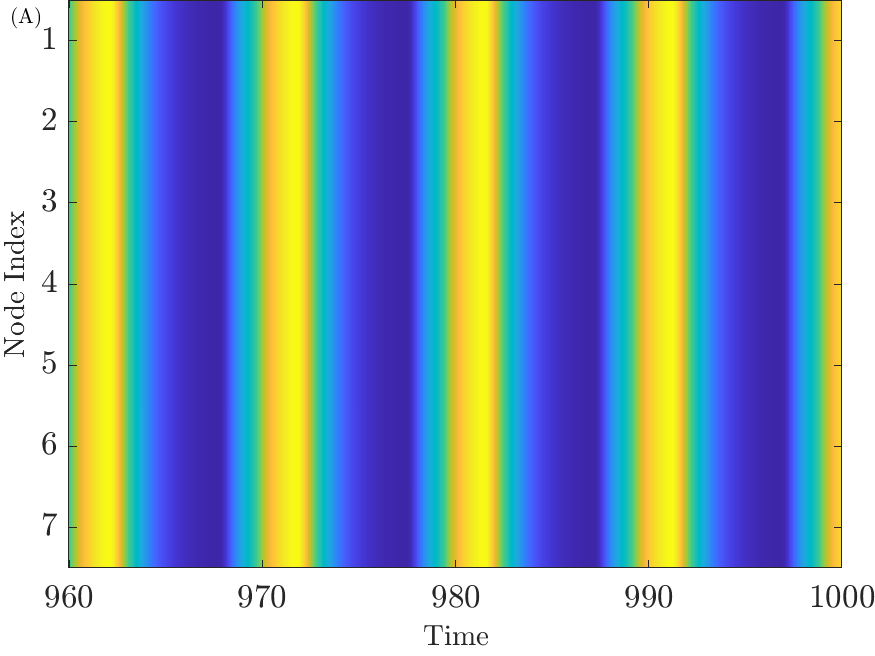}
\includegraphics[width=6cm]{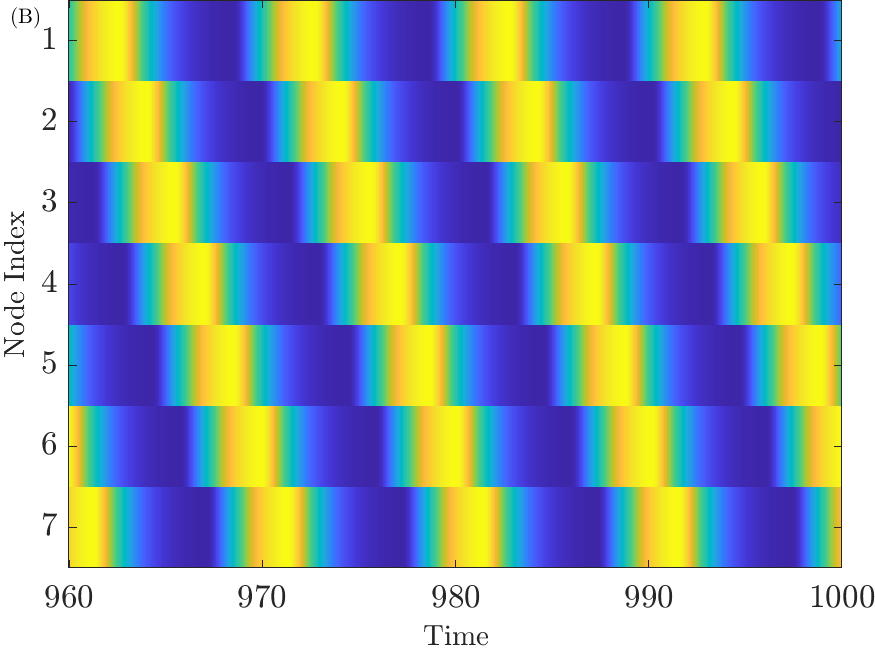}
\includegraphics[width=6cm]{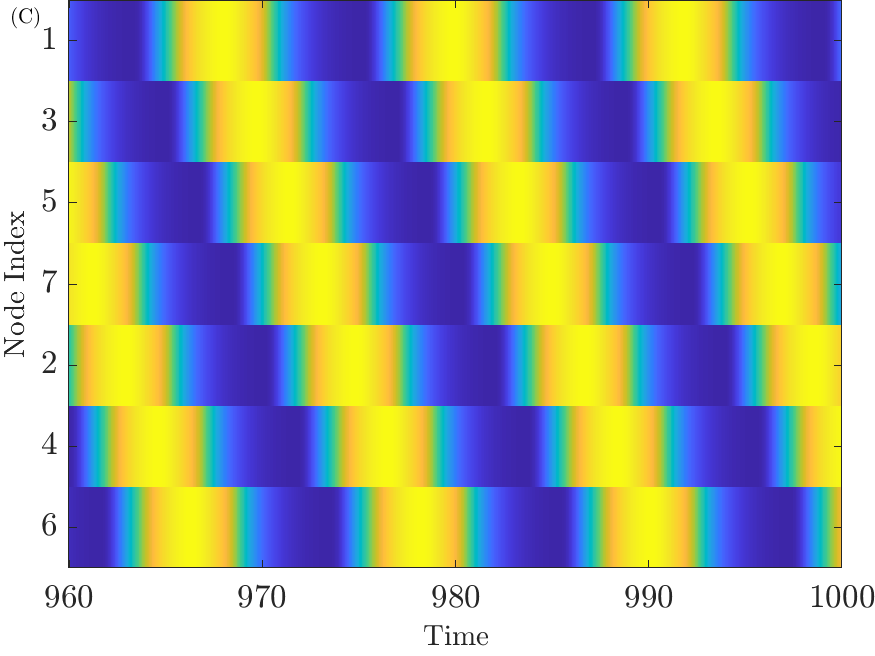}\\
\includegraphics[width=6cm]{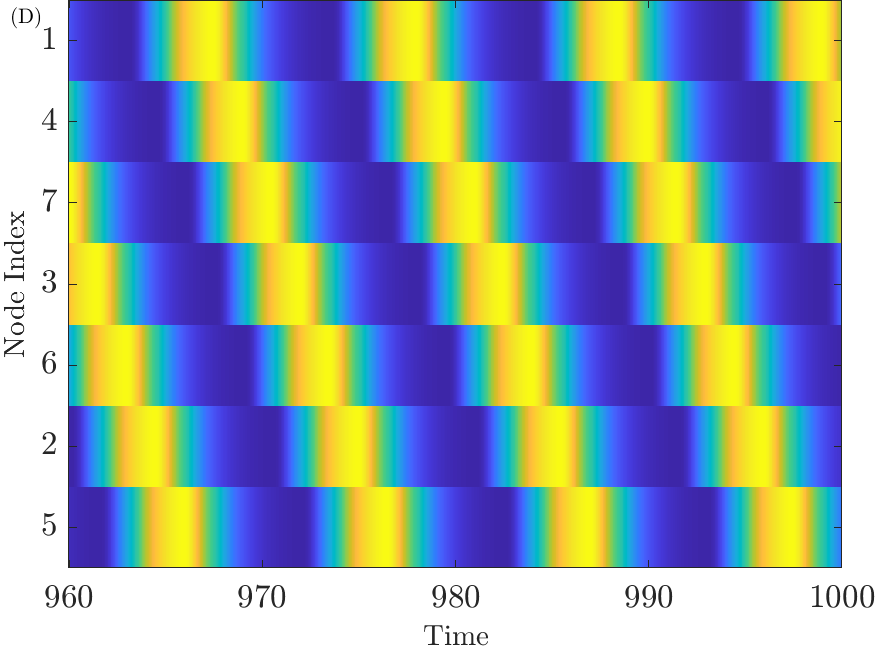}
\includegraphics[width=6cm]{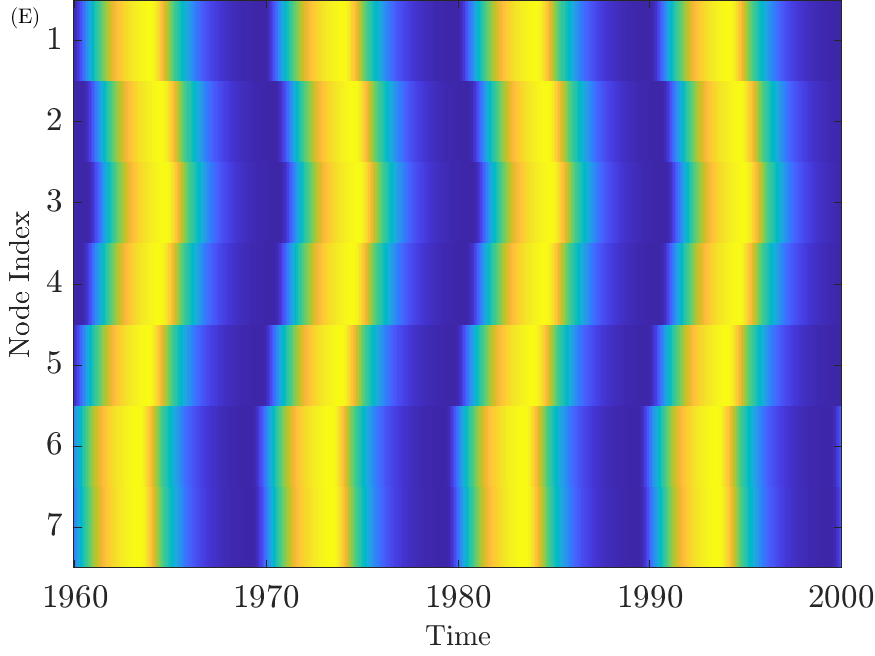}
\includegraphics[width=6cm]{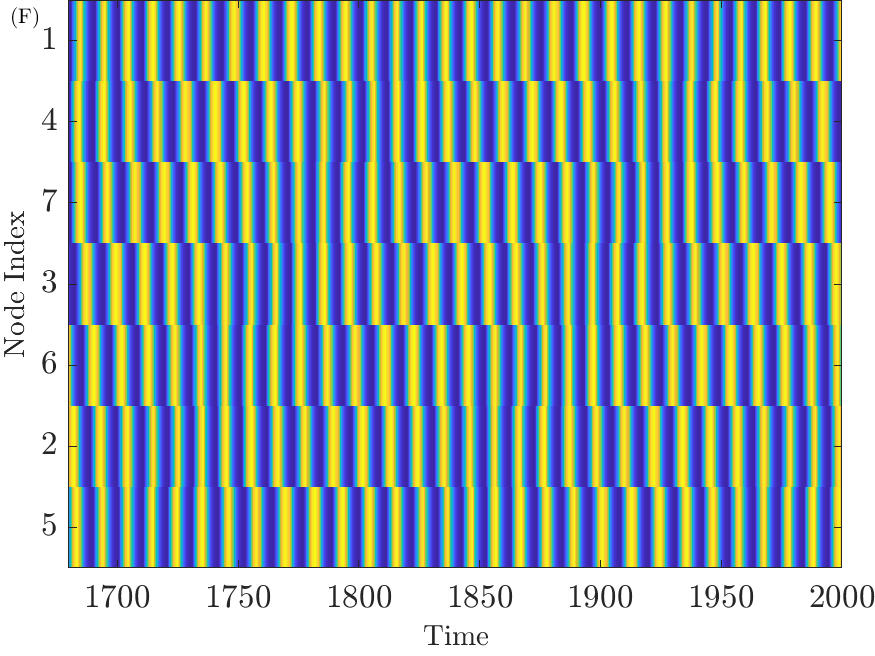}
\caption{Simulations of the Wilson--Cowan model (\ref{network}) with 7 nodes. (A) synchrony with $\tau_{\text{inter}}=5.3$,  (B) nearest-neighbour travelling wave with $\tau_{\text{inter}}=1.3$, (C) next-nearest neighbour travelling wave with $\tau_{\text{inter}}=2.3$, (D) next-next-nearest neighbour travelling wave with $\tau_{\text{inter}}=1.6$, (E) perturbed in-phase oscillation with $\tau_{\text{inter}}=6.2$, (F) Modulated travelling wave with $\tau_{\text{inter}}=6.5$. The node index has been shuffled to illustrate the various travelling waves. Here, $\tau_{\text{intra}}=1.5$ and other parameters fixed as in Fig.~\ref{Fig:WilsonCowan}.}
\label{fig:WCsim}
\end{figure}

Here we consider the Wilson--Cowan ring network given by (\ref{network}) with connectivity $\epsilon w_{ij}=0.2C\exp(-2 \operatorname{dist}(i,j))$, $w_{ii}=0$, and $C$ normalising the connectivity as in Section \ref{numerics_primer}. We first simulated the network with $\tau_{\text{inter}}=0.2$, see Fig.~\ref{fig:WCsim}A, and then determined the branch of synchronous solutions using numerical continuation. This synchronous branch exists throughout the domain we consider, i.e., from $\tau_{\text{inter}}=0.2$ to $\tau_{\text{inter}}=14.0$ and beyond. The profile varied only slightly, but in the plot of the period (not shown), we see it vary up and down, and then the branch starts to fold backwards near $\tau_{\text{inter}}=13.8$. Such folds may hint at bistability of the synchronous solution, but see below.

We then added small perturbations to the initial condition while varying $\tau_{\text{inter}}$, and travelling wave patterns and incoherent motion would emerge for certain values. For instance, we see travelling waves to nearest neighbour patterns ($q=1$), see Fig.~\ref{fig:WCsim}B or (next-)-next nearest nodes (corresponding to $q=2$ and $q=3$), see Figs. \ref{fig:WCsim}C,D. In between, there are also less coherent patterns such as a nearly synchronous mode and a torus-like solution, see Figs. \ref{fig:WCsim}E,F. The latter shows an additional phase of the oscillation travelling over the entire network.

\begin{figure}[ht!]
\includegraphics[width=16cm]{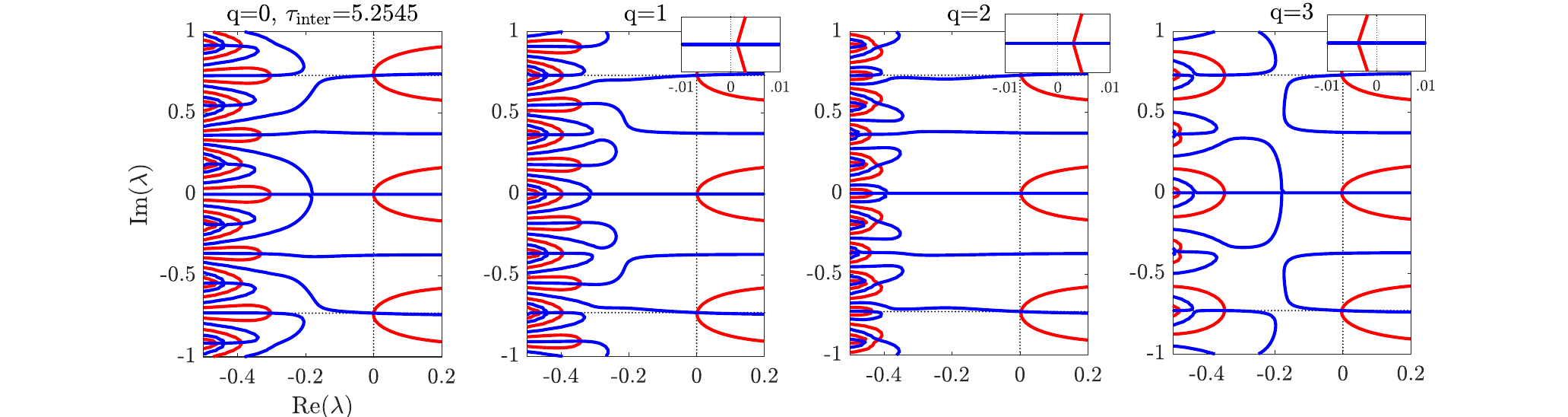}\\
\includegraphics[width=16cm]{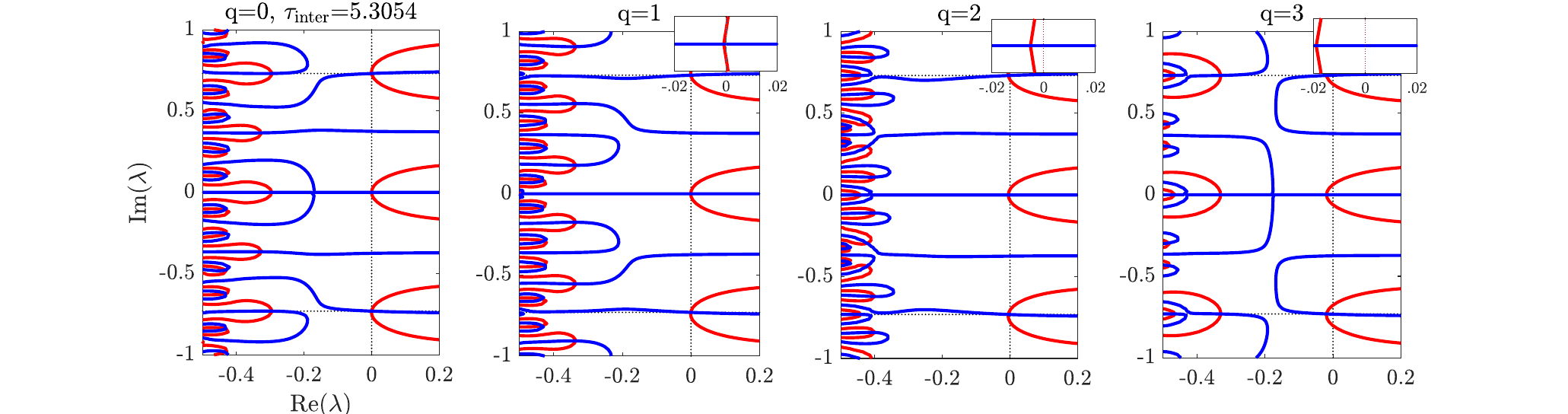}
\caption{For each $q$, the evaluation of the zero level sets  of $\mathcal{E}_q(\lambda)$ given by (\ref{Eq:RingStab}) for real (blue) and imaginary (red) parts. (Top) The insets show that the modes for $q=2$ and $q=3$ exhibit an instability for $\tau_{\text{inter}}=5.2545$; (Bottom) All eigenvalues are in the left half of the complex plane for $\tau_{\text{inter}}=5.3054$.}
\label{fig:WC_Spectrum_single}
\end{figure}

For the Wilson--Cowan ring network, each node has two state variables. Making use of the tensor approach developed in Sections \ref{WCsingledelay} and \ref{Sec:Ring},
we can compute the spectrum of the synchronous orbit using harmonic balance and the decomposition into modes. As for the primer, we show the spectrum before and after the instability. We selected $M=80$ Fourier modes because we found the spectrum was otherwise inaccurate. The cause is the sharper profile of the orbit that also needs more modes for an acceptable approximation, in contrast to the primer for which $M=30$ appeared sufficient for convergence. We select the instability near $\tau_{\text{inter}}=5.29$ as this point marks a region with physiologically plausible parameters, i.e. $\tau_{\text{inter}}\geq \tau_{\text{intra}}$. Evaluating (\ref{Eq:RingStab}) for each $q$, we observe a few real eigenvalues in the left half plane for $\tau_{\text{inter}}=5.3054$, and two unstable ones for $\tau_{\text{inter}}=5.2545$ for $q=2$ and $q=3$, see Figure \ref{fig:WC_Spectrum_single}. Simulations starting from the synchronous solution with a perturbation converge quickly to the next-nearest neighbour wave for $\tau_{\text{inter}}\leq 5.29$.

\begin{figure}[ht!]
\begin{center}
\includegraphics[width=7cm]{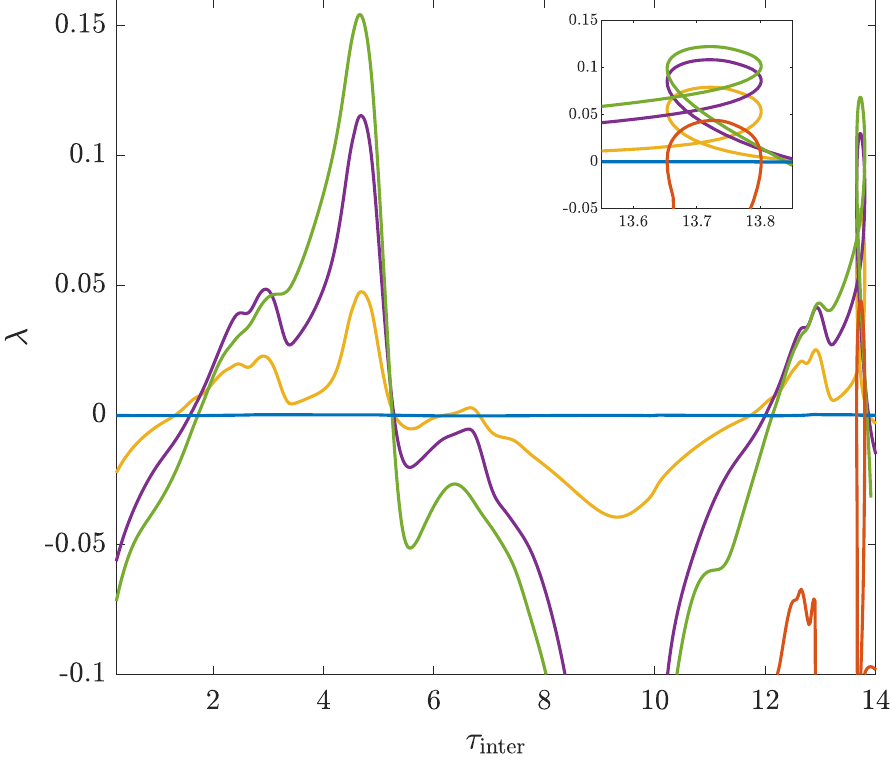}
\caption{For each $q$, the maximal eigenvalue as a function of $\tau_{\text{inter}}$ is shown. Line colour: blue (trivial, horizontal), nontrivial eigenvalue for $q=0$ (red), $q=1$ (yellow), $q=2$ (purple), $q=3$ (green). The inset shows that the branch makes a loop.} \label{fig:WC_Spectrum_full}
\end{center}
\end{figure}

Next, we looked at the maximal eigenvalues $\lambda_{q}$ along the entire branch, see Fig. \ref{fig:WC_Spectrum_full}. This $\lambda_{q}$ turned out to be real until $\tau_{\text{inter}}=13.91$, and we only used bisection to determine the largest $\lambda$. We see the first regions of (in)stability of the synchronous solution are demarcated by $\tau_{\text{inter}}\in\{1.3, 5.29, 6.12, 6,83, 11.72, 13.89\}$. In the region where the branch makes a loop, i.e. $13.65\leq \tau_{\text{inter}}\leq 13.80$, one may expect bistability of two different synchronous solutions within the synchronisation manifold $x_{i}(t)=x_{j}(t)$ for all $i,j=1, \ldots,N$. However, the other modes $ q=1,2,3$ are unstable, and so, there is no bistability of two synchronous solutions in the full system. 


\section{Discussion\label{Sec:Discussion}}

In this paper we have developed an applied mathematics approach for determining periodic orbits and their linear stability in neural mass networks with multiple delays. This is motivated by the growing number of computational studies of such networks in the context of large scale brain dynamics as seen in modern neuroimaging studies, although is equally applicable to other networks in the same universality class, such as physical models of lasers \cite{Soriano2013}, connected vehicle systems \cite{Szalai2013}, large ecosystems \cite{Pigani2022}, deep neural networks (that can be \textit{folded} into a single node with multiple time-delayed feedback loops) \cite{Stelzer2021}, and indeed other types of neural networks \cite{Cakan2014,Rahman2015}, including those for predictive coding \cite{Faye2023}, to name just a few.
Although the analytical and numerical study of DDEs has a strong history, this has focused mainly on \textit{node} as opposed to \textit{network} dynamics, and typically with a small number of delays. Here, we have shown that the harmonic balance approach combines nicely with Floquet theory to yield an analytic framework for obtaining results at the network level that can be developed into a practical numerical tool for network bifurcation analysis. This has allowed us to show how space-dependent delays in a Wilson--Cowan ring network of arbitrary size can destabilise the synchronous state in favour of travelling periodic waves, alternating anti-phase solutions, cluster states, and more exotic behaviours.
To complement the network Floquet theory that we have developed here, it would be interesting to construct the normal form for an emerging instability, say by using and extending ideas from the theory of \textit{Synergetics} \cite{Haken1993}. However, a more immediate challenge is to move away from ring structures of identical nodes and to consider delays that are truly heterogeneous as opposed to distance-dependent. Regarding the study of more general networks, it may be interesting to explore the use of techniques from computational group theory in revealing hidden symmetries that have allowed the study of phase-locked and cluster states in non-delayed networks \cite{Pecora2014}.
One might then develop a network lag operator, generalising the one introduced here, more apt for treating the combined space-time (strength \textit{and} delay) structure of human connectome data. On top of this, it is well to bear in mind that white matter is itself \textit{plastic}, see e.g., \cite{Fields2008,Gibson2014,Vivo2019} from an experimental perspective and e.g., \cite{Noori2020,Park2020,Talidou2022,Pajevic2023} from a modelling perspective. Such plasticity gives rise to state-dependent delays, and it seems fair to say that the mathematical study of state-dependent DDEs is in its relative infancy, see e.g., \cite{Humphries2022}. All of the above are topics of ongoing research and will be reported upon elsewhere.

\bibliographystyle{plain}
\bibliography{mnaHarmonicBalance}

\begin{thebibliography}{10}

\bibitem{Abeysuriya2018}
R~G Abeysuriya, J~Hadida, S~N Sotiropoulos, S~Jbabdi ansd R~Becker, B~A~E Hunt,
  M~J Brookes, and M~W Woolrich.
\newblock A biophysical model of dynamic balancing of excitation and inhibition
  in fast oscillatory large-scale networks.
\newblock {\em PLoS Computational Biology}, 14(2):e1006007, 2018.

\bibitem{Ahmadizadeh2016}
S~Ahmadizadeh, D~Ne\v{s}i\'c, D~R Freestone, and D~B Grayden.
\newblock On synchronization of networks of {Wilson-Cowan} oscillators with
  diffusive coupling.
\newblock {\em Automatica}, 71:169--178, 2016.

\bibitem{AlDarabsah2021}
I~Al-Darabsah, L~Chen, W~Nicola, and S~A Campbell.
\newblock The impact of small time delays on the onset of oscillations and
  synchrony in brain networks.
\newblock {\em Frontiers in Systems Neuroscience}, 15:688517, 2021.

\bibitem{Byrne2022}
\'A Byrne, J~Ross, R~Nicks, and S~Coombes.
\newblock {Mean-field models for EEG/MEG: from oscillations to waves}.
\newblock {\em Brain Topography}, 35:36--53, 2022.

\bibitem{Cabral2022}
J~Cabral, F~Castaldo, J~Vohryzek, V~Litvak, C~Bick, R~Lambiotte, K~Friston, M~L
  Kringelbach, and G~Deco.
\newblock Metastable oscillatory modes emerge from synchronization in the brain
  spacetime connectome.
\newblock {\em Nature Communications Physics}, 5:184, 2022.

\bibitem{Cakan2014}
C~Cakan, J~Lehnert, and E~Sch\"oll.
\newblock Heterogeneous delays in neural networks.
\newblock {\em The European Physical Journal B}, 87:54, 2014.

\bibitem{Campbell2007}
S~A Campbell.
\newblock {\em Time Delays in Neural Systems}, pages 65--90.
\newblock Springer Berlin Heidelberg, Berlin, Heidelberg, 2007.

\bibitem{Castaldo2023}
F~Castaldo, F~{P{\'a}scoa dos Santos}, R~C Timms, J~Cabral, J~Vohryzek, G~Deco,
  M~Woolrich, K~Friston, P~Verschure, and V~Litvak.
\newblock Multi-modal and multi-model interrogation of large-scale functional
  brain networks.
\newblock {\em NeuroImage}, 277:120236, 2023.

\bibitem{Conti2019}
F~Conti and R~A {Van Gorder}.
\newblock {The role of network structure and time delay in a metapopulation
  Wilson--Cowan model}.
\newblock {\em Journal of Theoretical Biology}, 477:1--13, 2019.

\bibitem{Coombes2018}
S~Coombes, Y~M Lai, M~{\c{S}}ayli, and R~Thul.
\newblock Networks of piecewise linear neural mass models.
\newblock {\em European Journal of Applied Mathematics}, 29:869--890, 2018.

\bibitem{Coombes2009}
S~Coombes and C~Laing.
\newblock Delays in activity-based neural networks.
\newblock {\em Philosophical Transactions of the Royal Society A},
  367:1117--1129, 2009.

\bibitem{Coombes2023}
S~Coombes and K~C~A Wedgwood.
\newblock {\em Neurodynamics: An Applied Mathematics Perspective}, volume~75 of
  {\em Texts in Applied Mathematics}.
\newblock Springer, 2023.

\bibitem{Vivo2019}
L~{de Vivo} and M~Bellesi.
\newblock The role of sleep and wakefulness in myelin plasticity.
\newblock {\em Glia}, 67:2142--2152, 2019.

\bibitem{Deco2009}
G~Deco, V~Jirsa, A~R McIntosh, O~Sporns, and R~K{\"o}tter.
\newblock {Key role of coupling, delay, and noise in resting brain
  fluctuations}.
\newblock {\em Proceedings of the National Academy of Sciences},
  106:10302--10307, 2009.

\bibitem{Detroux2015}
T~Detroux, L~Renson, L~Masset, and G~Kerschen.
\newblock The harmonic balance method for bifurcation analysis of large-scale
  nonlinear mechanical systems.
\newblock {\em Computer Methods in Applied Mechanics and Engineering},
  296:18--38, 2015.

\bibitem{MATCONT}
A~Dhooge, W~Govaerts, Y~A Kuznetsov, H~G~E Meijer, and B~Sautois.
\newblock {New features of the software MatCont for bifurcation analysis of
  dynamical systems}.
\newblock {\em Mathematical and Computer Modelling of Dynamical Systems},
  14:147--175, 2008.

\bibitem{Engelborghs2002}
K~Engelborghs, T~Luzyanina, and D~Roose.
\newblock {Numerical bifurcation analysis of delay differential equations using
  DDE-BIFTOOL}.
\newblock {\em ACM Transactions on Mathematical Software (TOMS)}, 28:1--21,
  2002.

\bibitem{Erneux2009}
T~Erneux.
\newblock {\em Applied Delay Differential Equations}.
\newblock Surveys and Tutorials in the Applied Mathematical Sciences. Springer,
  2009.

\bibitem{Faye2023}
G~Faye, G~Fouilh\'e, and R~{VanRullen}.
\newblock Mathematical derivation of wave propagation properties in
  hierarchical neural networks with predictive coding feedback dynamics.
\newblock {\em Bulletin of Mathematical Biology}, 85(80), 2023.

\bibitem{Fields2008}
R~D Fields.
\newblock White matter in learning, cognition and psychiatric disorders.
\newblock {\em Trends in Neurosciences}, 31:361--370, 2008.

\bibitem{Gibson2014}
E~M Gibson, D~Purger, C~W Mount, A~K Goldstein, G~L Lin, L~S Wood, I~Inema, S~E
  Miller, G~Bieri, J~B Zuchero, B~A Barres, P~J Woo, H~Vogel, and M~Monje.
\newblock Neuronal activity promotes oligodendrogenesis and adaptive
  myelination in the mammalian brain.
\newblock {\em Science}, 344:1252304, 2014.

\bibitem{Gross2019}
J~Gross.
\newblock Magnetoencephalography in cognitive neuroscience: A primer.
\newblock {\em Neuron}, 104:189--204, 2019.

\bibitem{Haken1993}
H~Haken.
\newblock {\em Advanced Synergetics}.
\newblock Springer-Verlag, New York, 1993.

\bibitem{Hoppensteadt97}
F~C Hoppensteadt and E~M Izhikevich.
\newblock {\em Weakly Connected Neural Networks}.
\newblock Springer-Verlag, Heidelberg, Germany, 1997.

\bibitem{Humphries2022}
A~R Humphries, B~Krauskopf, S~Ruschel, and J~Sieber.
\newblock Nonlinear effects of instantaneous and delayed state dependence in a
  delayed feedback loop.
\newblock {\em Discrete and Continuous Dynamical Systems - B}, 27:7245--7273,
  2022.

\bibitem{Kotani2012}
K~Kotani, I~Yamaguchi, Y~Ogawa, Y~Jimbo, H~Nakao, and G~B Ermentrout.
\newblock Adjoint method provides phase response functions for delay-induced
  oscillations.
\newblock {\em Physical Review Letters}, 109:044101, Jul 2012.

\bibitem{Krylov1943}
N~M Krylov and N~N Bogoliubov.
\newblock {\em Introduction to Non-linear Mechanics}.
\newblock Princeton University Press, 1943.

\bibitem{Kutchko2013}
K~Kutchko and F~Fr\"ohlich.
\newblock Emergence of metastable state dynamics in interconnected cortical
  networks with propagation delays.
\newblock {\em PLoS Computational Biology}, 9:e1003304, 2013.

\bibitem{Liu2010}
L~Liu and T~Kalm\'ar-Nagy.
\newblock High-dimensional harmonic balance analysis for second-order
  delay-differential equations.
\newblock {\em Journal of Vibration and Control}, 16:1189--1208, 2010.

\bibitem{MacDonald1995}
N~MacDonald.
\newblock Harmonic balance in delay-differential equations.
\newblock {\em Journal of Sound and Vibration}, 186:649--656, 1995.

\bibitem{MeijerDercoleOldeman2009}
H.G.E. Meijer, F.~Dercole, and B.~Oldeman.
\newblock {\em Numerical Bifurcation Analysis}, pages 6329--6352.
\newblock Springer New York, New York, NY, 2009.

\bibitem{Nakagawa2023}
T~T Nakagawa, M~Woolrich, H~Luckhoo, M~Joensson, H~Mohseni, M~L Kringelbach,
  V~Jirsa, and G~Deco.
\newblock How delays matter in an oscillatory whole-brain spiking-neuron
  network model for {MEG} alpha-rhythms at rest.
\newblock {\em NeuroImage}, 87:383--394, 2014.

\bibitem{Nicks2024}
R~Nicks, R~Allen, and S~Coombes.
\newblock Phase and amplitude responses for delay equations using harmonic
  balance.
\newblock {\em Physical Review E}, 110:L012202, 2024.

\bibitem{Noori2020}
R~Noori, D~Park, J~D Griffiths, S~Bells, P~W Frankland, D~Mabbott, and
  J~Lefebvre.
\newblock {Activity-dependent myelination: A glial mechanism of oscillatory
  self-organization in large-scale brain networks}.
\newblock {\em Proceedings of the National Academy of Sciences USA},
  117:13227--13237, 2020.

\bibitem{Novicenko2012}
V~Novi\v{c}enko and K~Pyragas.
\newblock Phase reduction of weakly perturbed limit cycle oscillations in
  time-delay systems.
\newblock {\em Physica D: Nonlinear Phenomena}, 241:1090--1098, 2012.

\bibitem{Otto2018}
A~Otto, G~Radons, D~Bachrathy, and G~Orosz.
\newblock Synchronization in networks with heterogeneous coupling delays.
\newblock {\em Physical Review E}, 97:012311, 2018.

\bibitem{Pajevic2023}
S~Pajevic, D~Plenz, P~J Basser, and R~D Fields.
\newblock Oligodendrocyte-mediated myelin plasticity and its role in neural
  synchronization.
\newblock {\em eLife}, 12:e81982, 2023.

\bibitem{Park2020}
S~H Park and J~Lefebvre.
\newblock {Synchronization and resilience in the Kuramoto white matter network
  model with adaptive state-dependent delays}.
\newblock {\em The Journal of Mathematical Neuroscience}, 10(1):16, 2020.

\bibitem{Pecora2014}
L~M Pecora, F~Sorrentino, A~M Hagerstrom, T~E Murphy, and R~Roy.
\newblock Cluster synchronization and isolated desynchronization in complex
  networks with symmetries.
\newblock {\em Nature Communications}, 5, 2014.

\bibitem{Petkoski2019}
S~Petkoski and V~K Jirsa.
\newblock Transmission time delays organize the brain network synchronization.
\newblock {\em Philosophical Transactions of the Royal Society A}, 37:20180132,
  2019.

\bibitem{Petkoski2016}
S~Petkoski, A~Spiegler, T~Proix, P~Aram, J-J Temprado, and V~K Jirsa.
\newblock Heterogeneity of time delays determines synchronization of coupled
  oscillators.
\newblock {\em Physical Review E}, 94:012209, Jul 2016.

\bibitem{Pigani2022}
E~Pigani, D~Sgarbossa, S~Suweis, A~Maritan, and S~Azaele.
\newblock Delay effects on the stability of large ecosystems.
\newblock {\em Proceedings of the National Academy of Sciences USA},
  119:e2211449119, 2022.

\bibitem{Pinder2023}
I~Pinder, M~R Nelson, and J~J Crofts.
\newblock {Bifurcations and synchrony in a ring of delayed Wilson--Cowan
  oscillators}.
\newblock {\em Proceedings of the Royal Society A}, 470:20230313, 2023.

\bibitem{Rackauckas2017}
C~Rackauckas and Q~Nie.
\newblock {DifferentialEquations.jl--A performant and feature-rich ecosystem
  for solving differential equations in Julia}.
\newblock {\em Journal of Open Research Software}, 5(1):15, 2017.

\bibitem{Rahman2015}
B~Rahman, K~B Blyuss, and Y~N Kyrychko.
\newblock Dynamics of neural systems with discrete and distributed time delays.
\newblock {\em SIAM Journal on Applied Dynamical Systems}, 14:2069--2095, 2015.

\bibitem{Ruschel2025}
S~Ruschel and A~Giraldo.
\newblock Master stability for traveling waves on networks.
\newblock {\em Physical Review Letters}, 134:257201, 2025.

\bibitem{Sayli2024}
M~{\c{S}ayli} and S~Coombes.
\newblock Understanding the effect of white matter delays on large scale brain
  synchrony.
\newblock {\em Communications in Nonlinear Science and Numerical Simulation},
  131:107803, 2024.

\bibitem{Simmendinger1999}
C~Simmendinger, A~Wunderlin, and A~Pelster.
\newblock Analytical approach for the {F}loquet theory of delay differential
  equations.
\newblock {\em Physical Review E}, 59:5344--5353, 1999.

\bibitem{Soriano2013}
M~C Soriano, J~Garc\'{\i}a-Ojalvo, C~R Mirasso, and I~Fischer.
\newblock {Complex photonics: Dynamics and applications of delay-coupled
  semiconductors lasers}.
\newblock {\em Reviews of Modern Physics}, 85:421--470, 2013.

\bibitem{Stelzer2021}
F~Stelzer, A~R\"ohm, R~Vicente, I~Fischer, and S~Yanchuk.
\newblock Deep neural networks using a single neuron: folded-in-time
  architecture using feedback-modulated delay loops.
\newblock {\em Nature Communications}, 12:5164, 2021.

\bibitem{Sun2023}
P~Sun, X~Zhao, X~Yu, Q~Huang, Z~Feng, and J~Zhou.
\newblock Incremental harmonic balance method for multi-harmonic solution of
  high-dimensional delay differential equations: Application to
  crossflow-induced nonlinear vibration of steam generator tubes.
\newblock {\em Applied Mathematical Modelling}, 118:818--831, 2023.

\bibitem{Szalai2013}
R~Szalai and G~Orosz.
\newblock Decomposing the dynamics of heterogeneous delayed networks with
  applications to connected vehicle systems.
\newblock {\em Physical Review E}, 88:040902, 2013.

\bibitem{Talidou2022}
A~Talidou, P~W Frankland, D~Mabbott, and J~Lefebvre.
\newblock {Homeostatic coordination and up-regulation of neural activity by
  activity-dependent myelination}.
\newblock {\em Nature Computational Science}, 2:665--676, 2022.

\bibitem{Tewarie2019}
P~Tewarie, R~Abeysuriya, \'A Byrne, G~C O'Neill, S~N Sotiropoulos, M~J Brookes,
  and S~Coombes.
\newblock How do spatially distinct frequency specific {MEG} networks emerge
  from one underlying structural connectome? {T}he role of the structural
  eigenmodes.
\newblock {\em NeuroImage}, 186:211--220, 2019.

\bibitem{Ton2014}
R~Ton, G~Deco, and A~Daffertshofer.
\newblock {Structure-Function Discrepancy: Inhomogeneity and Delays in
  Synchronized Neural Networks}.
\newblock {\em PLoS Computational Biology}, 10:e1003736, 2014.

\bibitem{Wilson72}
H~R Wilson and J~D Cowan.
\newblock Excitatory and inhibitory interactions in localized populations of
  model neurons.
\newblock {\em Biophysical Journal}, 12:1--24, 1972.

\end{thebibliography}

\end{document}